\documentclass{amsart}
\usepackage{amssymb}
\usepackage{amscd}
\usepackage{amsmath}

\theoremstyle{plain}
\newtheorem{theorem}{Theorem}[section]
\newtheorem{proposition}{Proposition}[section]
\newtheorem{lemma}{Lemma}[section]
\newtheorem{cor}{Corollary}[section]

\theoremstyle{definition}

\theoremstyle{remark}

\newtheorem{remark}{Remark}[section]

\begin{document}


\title[Moment Conditions and Support Theorems]%
{Moment Conditions and Support Theorems for Radon Transforms on Affine Grassmann Manifolds}

\author{Fulton B. Gonzalez}
\address{Department of Mathematics\\
Tufts University}
\email{fulton.gonzalez@tufts.edu}
\author{Tomoyuki Kakehi}
\address{Institute of Mathematics\\
 University of Tsukuba}
\email{kakehi@math.tsukuba.ac.jp}

\date{March 3, 2004}
\maketitle

\def\rar{\rightarrow}
\def\fk{\mathfrak}
\def\e{\frak e}
\def\g{\frak g}
\def\k{\frak k}
\def\a{\frak a}
\def\m{\frak m}
\def\n{\frak n}
\def\p{\frak p}
\def\u{\frak u}
\def\t{\frak t}
\def\h{\frak h}
\def\z{\frak z}
\def\d{\frak d}

\def\Pf{\mathrm{Pf}}
\def\rank{\mathrm{rank}}
\def\Pr{\mathrm{Pr}}
\def\Vol{\mathrm{Vol}}

\numberwithin{equation}{section}


\begin{abstract}
Let $G(p,n)$ and $G(q,n)$ be the affine Grassmann manifolds
of $p$- and $q$- planes in ${\mathbb R}^n$, respectively,
and let $\mathcal{R}^{\, (p,q)}$
be the Radon transform from smooth functions on $G(p,n)$
to smooth functions on $G(q,n)$ arising from the inclusion incidence relation.
When $p<q$ and $\dim G(p,n) = \dim G(p,n)$, 
we present a range characterization theorem
for $\mathcal{R}^{\, (p,q)}$ via moment conditions.
We then use this range result to prove a support theorem
for $\mathcal{R}^{\, (p,q)}$.
This complements a previous range characterization theorem
for $\mathcal{R}^{\, (p,q)}$ via differential equations
when $\dim G(p,n) < \dim G(p,n)$.
We also present a support theorem in this latter case.
\end{abstract}

\begin{section}{Introduction.}

In this paper, we present a range characterization of
Radon transforms on affine Grassmann manifolds via moment conditions.
Our objective is to generalize the moment conditions for 
the classical Radon transform on $\mathbb{R}^n$ to affine Grassmannians;
our results complement a previously obtained range characterization
for these transforms using invariant differential equations.

Specifically, we consider the transform $\mathcal{R}^{\, (p,q)}$ from
smooth functions on the space $G(p,n)$ of $p$-planes in $\mathbb{R}^n$
to smooth functions on the space $G(q,n)$ of $q$-planes in $\mathbb{R}^n$
arising from the inclusion incidence relation.
Throughout this paper, we assume that $p<q$.
A $p$-plane $\ell$ and a $q$-plane $\xi$ are incident if $\ell \subset \xi$;
the Radon transform is explicitly defined by
\begin{equation}
\mathcal{R}^{\, (p,q)} \, f (\xi) \; = \;
   \int_{\widehat{\xi}} \; f(\ell) \, d_{\xi} \ell.
\end{equation}
when $\xi \in G(q,n)$ and $f$ is a function on $G(p,n)$.
Here $d_{\xi} \ell$ is a canonical measure on the set
$\widehat{\xi} = \{ \; \ell \in G(p,n) \; | \ell \; \text{is incident to} \; \xi \; \}$
invariant under all Euclidean motions preserving $\eta$.
In particular, when $p=0$, $\mathcal{R}^{\, (0,q)}$ reduces to 
the classical $q$-plane transform on $\mathbb{R}^n$.
We will be working in the category $\mathcal{S} (G(p,n))$ of Schwartz class functions
on $G(p,n)$, as defined in \cite{Richter}.
The following theorem [GK 1, Theorem 7.7] gives a range characterization of
$\mathcal{R}^{\, (p,q)}$ in the case when $\dim G(p,n) < \dim G(q,n)$:

\begin{theorem}\label{T:diffrangechar}
Assume that $p < q$ and $\dim G(p,n) < \dim G(q,n)$.
Then there exists a differential operator $\Omega$ of order $2p+4$,
invariant under the Euclidean motion group $M(n)$, such that for any
$\varphi \in \mathcal{S} (G(q,n))$,
\begin{equation}
\varphi \in \mathcal{R}^{\, (p,q)} \mathcal{S} (G(p,n))
\; \Longleftrightarrow \; \Omega \varphi =0.
\end{equation}
\end{theorem}

This theorem generalizes the range characterization of the $k$-plane transform
on $\mathbb{R}^n$, when $k < n-1$, via a second order ultrahyperbolic system
(\cite{Richter}, \cite{Gonzalez-TAMS}) a single 4th order $M(n)$-invariant differential operator
(\cite{Gonzalez-Annalen}).

The classical Radon transform, on the other hand, has its range given by moment conditions,
as specified in \cite{Gelfand-Graev-Vilenkin}, \cite{Ludwig}, \cite{Helgason-1965}.
(See the last reference for a complete proof.)
The Grassmannian analogue to this would correspond to the transform $\mathcal{R}^{\, (p,q)}$
with $\dim G(p,n) = \dim G(q,n)$.
Since the dimensions coincide, we would expect the range to be
specified by appropriate moment conditions
as well.
This is the content of our main result, Theorem~\ref{T:momentconditions} below.
As with the proof in \cite{Helgason-1965}, the crucial element of our proof consists of
justifying the smoothness of a certain ``partial Fourier transform" on the space
$G(p,n)$.
(Propositions~\ref{T:f-smoothness} and \ref{T:f-rapid-decrease} below.)

Now Helgason's geometric proof of the support theorem for the
classical Radon transform is well-known (\cite{Helgason-1965}; less
well-known is the fact that it is also a consequence of the forward
(easy) moment conditions and a polar-coordinate version of the
Paley-Wiener Theorem (\cite{GGA}; Theorem~\ref{T:Paleypolar1} below)). 
We give an extension of this theorem to affine Grassmannians (Theorem~\ref{T:Paleypolar2})
and use it to prove a support theorem
(Theorem~\ref{T:Grassmanniansupport}) for the transform $\mathcal R^{(p,q)}$
in the case when $\dim (G(p,n))\,=\,\dim (G(q,n))$.  The injectivity of $\mathcal R^{(p,q)}$
is then used to prove the support theorem
(Theorem~\ref{th:support-theorem-1}) when  $\dim (G(p,n))\,<\,\dim (G(q,n))$.

\end{section}


\begin{section}{Moment Conditions and the Support Theorem 
for the Classical Radon Transform, Revisited.}

To clarify our exposition, it will be instructive to briefly summarize Helgason's proof 
of the range characterization
of the classical Radon transform by moment conditions.
For this, let $\mathcal{R} = \mathcal{R}^{\,
(0,n-1)}:f\,\mapsto\widehat f$ denote the classical Radon
transform, which maps functions in $\mathcal{S} (\mathbb{R}^n)$ to functions 
$\varphi \in \mathcal{S} (\mathbf{S}^{n-1} \times \mathbb{R})$
which satisfy $\varphi (\omega, p) = \varphi (-\omega, -p)$.
(See \cite{GGA}, pg.99 for the definition of $\mathcal{S} (\mathbf{S}^{n-1} \times \mathbb{R})$.)
We call such functions even and define
$\mathcal{S}_H (\mathbf{S}^{n-1} \times \mathbb{R})$ to be the vector space
of all even functions in $\mathcal{S} (\mathbf{S}^{n-1} \times \mathbb{R})$ satisfying
the following condition:
\begin{equation}
\begin{split}
\; &
\text{\it For any $k \in \mathbb{Z}^{+}$, there exists a homogeneous polynomial
$P_k$ of degree $k$ on $\mathbb{R}^n$ such that}\\
\; & \qquad\qquad \int_{-\infty}^{\infty} \; \varphi (\omega, p) \, p^k \; dp
\; = \; P_k (\omega),  \qquad \text{for} \; \; \omega \in \mathbf{S}^{n-1}.
\end{split}
\tag{H}\label{E:H}
\end{equation}

\begin{theorem}
\; (\cite{Helgason-1965}) \; \; 
$\mathcal{R} \mathcal{S} (\mathbb{R}^n) = \mathcal{S}_H (\mathbf{S}^{n-1} \times \mathbb{R})$.
\end{theorem}

It is an easy calculation to show that
$$
\int_{-\infty}^{\infty} \; \mathcal{R} f (\omega, p) \, p^k \; dp
\; = \;
\int_{\mathbb{R}^n} \;  f (x) \, \langle x, \omega \rangle^k \; dx.
$$
This shows that
$\mathcal{R} \mathcal{S} (\mathbb{R}^n) \subset
 \mathcal{S}_H (\mathbf{S}^{n-1} \times \mathbb{R})$.
Conversely, suppose that
$\varphi \in \mathcal{S}_H (\mathbf{S}^{n-1} \times \mathbb{R})$.
Define the ``partial Fourier transform'' $\widetilde{\varphi}$ on
$\mathbf{S}^{n-1} \times \mathbb{R}$,
\begin{equation}
\widetilde{\varphi} (\omega, s)
\; = \;
\int_{-\infty}^{\infty} \;  \varphi (\omega, p) \, e^{-ips} \; dp,
\qquad s \in \mathbb{R}, \; \omega \in \mathbf{S}^{n-1}.\label{E:PartialFT1}
\end{equation}
It is straightforward to prove that $\widetilde{\varphi}$ is an even function
in $\mathcal{S} (\mathbf{S}^{n-1} \times \mathbb{R})$.
In addition, the condition (\ref{E:H}) for $k=0$ shows that 
$\omega\mapsto\widetilde{\varphi} (\omega, 0)$ is constant so that there exists a unique function
$F$ on $\mathbb{R}^n$ for which
$F(s \omega ) = \widetilde{\varphi} (\omega, s)$.
Now the map $(\omega, s) \mapsto s \omega$ is a local diffeomorphism of
$\mathbf{S}^{n-1} \times (\mathbb{R} \setminus \{ 0 \})$ onto
$\mathbb{R}^n \setminus \{ 0 \}$, which shows that $F$ is smooth outside the origin.
To prove the smoothness of $F$ at the origin, it suffices to show that
the partial derivatives of $F$ are bounded on, say, the punctured unit ball
$B' = \{ \; x \in \mathbb{R}^n \; | \; 0 < || x || <1 \; \}$.
Fix $\epsilon >0$.
Now on the subset
$A_{n,\epsilon}\,=\,\{s \omega \in B' \, | \, 0 < s < 1, \; 
\omega = (\omega_1, \cdots, \omega_n) \in \mathbf{S}^{n-1},
\; \omega_n > \epsilon > 0 \}$
we can use $(s, \omega_1, \cdots, \omega_{n-1})$ as local coordinates:
repeated application of the chain rule shows that
\begin{equation}
\frac{\partial^k \; F}{\partial x_{i_1} \cdots \partial x_{i_k}} (x)
\; = \;
\sum_{j; i_1, \cdots, i_m} \; 
\frac{A_{j; i_1, \cdots, i_m} ( \omega_1, \cdots, \omega_{n-1})}{s^{k-j}} \,
\frac{\partial^m }{\partial \omega_{i_1} \cdots \partial \omega_{i_m}} \, 
\frac{\partial^j}{\partial s^j} \, F(s \omega )\label{E:polardiff}
\end{equation}
The coefficients $A_{j; i_1, \cdots, i_m} ( \omega_1, \cdots, \omega_{n-1})$
are smooth bounded functions of $( \omega_1, \cdots, \omega_{n-1})$ and the right hand
sum ranges over all $j \leq k$ and over all sequences
$i_1, \cdots, i_m$ in $\{ 1, 2, \cdots, n-1 \}$ where $m \leq k$.

We write $e^{-ips}\,=\,\sum_{l=0}^{k-1} (-ips)^l/l!\,+\,e_k(-ips)$
and apply condition (\ref{E:H}) to (\ref{E:PartialFT1}) to obtain
$$
F(s\omega)\;=\;\sum_{l=0}^{k-1}\frac{(-i)^l}{l!}P_l(s\omega)\,+\,\int_{-\infty}^\infty
\varphi(\omega,p)e_k(-ips)\,dp.
$$
Since the $P_l$ are polynomials of degree $l$, (\ref{E:polardiff})
implies that
\begin{equation}
\frac{\partial^k F}{\partial x_{i_1}\cdots\partial
x_{i_k}}(x)\;=\;\sum_{j;i_1,\ldots,i_m}
A_{j;i_1,\ldots,i_m}(\omega)\,\,\int_{-\infty}^\infty
\frac{\partial^m\varphi}{\partial\omega_{i_1}\cdots\partial\omega_{i_m}}(\omega,p)\,(-ip)^k\,
\frac{e_{k-j}(-ips)}{(-ips)^{k-j}}\,dp \label{E:polardiff2}
\end{equation}

Now $\varphi\in\mathcal{S}(S^{n-1}\times\mathbb R)$ and
$e_{k-j}(-it)/(-it)^{k-j}$ is bounded for all real $t$, so
(\ref{E:polardiff2}) shows
that the $k$th order derivatives of $F$ are bounded on the set
$A_{n,\epsilon}$.  From this, we deduce that the $k$th order
derivatives of $F$ are bounded on $B'$.  Hence $F\in C^\infty(\mathbb
R^n)$.  A routine calculation using, say, (\ref{E:polardiff}), shows
that $F\in\mathcal S(\mathbb R^n)$.   Denoting the
Fourier transform on $\mathbb R^n$ by  $\mathcal F$, let $f$ be the inverse Fourier
transform of $F$:  $F=\mathcal F(f)$.  For any $\omega,\,s$  the projection-slice theorem
says that $\mathcal Ff(s\omega)\,=\,\int_{-\infty}^\infty \widehat
f(\omega,p)\,e^{-ips}\,dp$.  The
injectivity of the Fourier transform on $\mathbb R$ then implies that
$\varphi\,=\,\widehat f$.\;\;\qedsymbol
\medskip

We now state the support theorem for the classical Radon transform in
the following form.

\begin{theorem}\label{T:classicalsupport} (\cite{Helgason-1965})
Let $f\in\mathcal S(\mathbb R^n)$
and suppose that $R>0$.  If $\widehat f(\omega,p)\,=\,0$ whenever
$|p|>R$, then $f$ has support in the closed ball $\bar B_R\,=\,\{x\in \mathbb
R^n\,|\;\|x\|\leq R\}$.
\end{theorem}

While the support theorem can be proved geometrically, it is also a
consequence of the forward  moment conditions $(H)$.  The key is
the following polar coordinate version of the Paley-Wiener theorem.  (See
exercise B1, Ch. 1 in \cite{GGA}.\footnote{The authors would like to thank Prof. S. 
Helgason for pointing out this exercise.})

\begin{theorem}(\cite{Helgason-1976})\label{T:Paleypolar1}
Let $R>0$.  The Fourier transform $f\mapsto
\mathcal Ff$ maps $\mathcal D(\bar B_R)$ onto the set of functions
$\mathcal Ff(\lambda \omega)=\psi(\omega,\lambda)\in C^\infty(S^{n-1}\times \mathbb R)$
satisfying the following conditions:
\item{(i)} For each $\omega\in S^{n-1}$, the function $\lambda\mapsto
\psi(\omega,\lambda)$ extends to a holomorphic function on $\mathbb C$
with the property that 
\medskip
$$
\sup_{\omega,\lambda}
\left|\psi(\omega,\lambda)\,(1+|\lambda|)^N\,e^{-R|\text{Im}\lambda|}\right|\;<\;\infty
\qquad\qquad\quad
N\in\mathbb Z^+.
$$
\item{(ii)} For each $k\in\mathbb Z^+$ and each homogeneous degree $k$
spherical harmonic function $h$ on $S^{n-1}$, the function
$$
\lambda\,\mapsto\,\lambda^{-k}\,\int_{S^{n-1}}\psi(\omega,\lambda)\,h(\omega)\,d\omega
$$
is even and holomorphic on $\mathbb C$ ($d\omega$ denoting area
measure on $S^{n-1}$).
\end{theorem}

To see how the support theorem follows from this, we take $f\in
\mathcal S(\mathbb R^n)$ satisfying $\widehat f(\omega,p)\,=\,0$ for
all $|p|>R$.  Then $\psi(\omega,s)\,=\,\mathcal F f(s\omega)$
satisfies condition (i) above by the easy part of the Paley-Wiener
theorem and the projection-slice theorem.  Now let $h$ be 
a homogeneous degree $k$ spherical harmonic
on $S^{n-1}$.  Then
$$
\begin{aligned}
\int_{S^{n-1}}\psi(\omega,s)\,h(\omega)\,d\omega\;&=\;\int_{S^{n-1}}\int_{-\infty}^\infty
\widehat f(\omega,p)\,e^{-ips}\,dp\,h(\omega)\,d\omega\\
&=\;\sum_{l=0}^{k-1}\frac
1{l!}\,\int_{S^{n-1}}\int_{-\infty}^\infty\widehat 
f(\omega,p)\,(-ips)^l\,dp\,h(\omega)\,d\omega \\
& \qquad\qquad\qquad +\;\int_{S^{n-1}}\int_{-\infty}^\infty\widehat
f(\omega,p)\,h(\omega)\,d\omega\,e_k(-ips)\,dp\\
&=\;\sum_{l=0}^{k-1}\frac
{(-is)^l}{l!}\,\int_{S^{n-1}}P_l(\omega)\,h(\omega)\,d\omega\\
&\qquad\qquad\qquad +\;\int_{S^{n-1}}\int_{-\infty}^\infty\widehat
f(\omega,p)\,h(\omega)\,d\omega\,e_k(-ips)\,dp
\end{aligned}
$$
by the forward moment conditions (H) for $\widehat f$.  Since
$P_l(\omega)$ is a sum of spherical harmonics of degree $\leq l$, the
sum on the right-hand side vanishes, so that
$$
s^{-k}\,\int_{S^{n-1}}\psi(\omega,s)\,h(\omega)\,d\omega\;=\;\int_{-\infty}^\infty
s^{-k}\,e_k(-ips)\,\int_{S^{n-1}} \widehat f(\omega,p)\,h(\omega)\,d\omega\,dp.
$$
The inner integral on the right is a smooth compactly supported
function of $p\in \mathbb R$ and $s\,\mapsto\,s^{-k}e_k(-ips)$ extends
to a holomorphic function on $\mathbb C$.  In addition, the right hand
side above is easily seen to be even in $s$.  Hence $\psi (\omega,s)$
satisfies condition (ii) in Theorem~\ref{T:Paleypolar1}, and so $f$ is
supported in the closed ball $\bar B_R$.
\end{section}


\begin{section}{The range of the Radon transform on affine
Grassmannians: the equal rank case}

We adopt the notation of \cite{Gonzalez-Kakehi-1} in what follows.
Let $G_{p,n}$ be the (compact) Grassmann manifold of $p$-dimensional
subspaces of $\mathbb R^n$.  Then $G_{p,n}=O(n)/K_p$, where
$K_p=O(p)\times O(n-p)$ is the subgroup of $O(n)$ fixing the $p$-plane
$\sigma_0\,=\,\mathbb R e_1\,+\cdots\,\mathbb R e_p$.

We assume that the Haar measure on all compact Lie groups, and the
invariant measures on their homogeneous spaces (with the exception of
the unit spheres), are normalized.   Let
$R_{p,q,n}:\,C^\infty(G_{p,n})\,\rightarrow\,C^\infty (G_{q,n})$ be the
Radon transform corresponding to the inclusion incidence relation
between $p$\,- \,and $q$\,-\,dimensional subspaces of $\mathbb R^n$.
Then $R_{p,q,n}$ is a linear bijection when $\rank (G_{p,n})\,=\,\rank
(G_{q,n})$ (\cite{Grinberg}); when $\rank (G_{p,n})\,<\,\rank
(G_{q,n})$, $R_{p,q,n}$ is injective and the range $R_{p,q,n}\,C^\infty(G_{p,n})$ is the subspace of
$C^\infty (G_{q,n})$ annihilated by an $O(n)$-invariant differential
operator of order $2\, \rank (G_{p,n})\,+\,2$.\quad (\cite{Kakehi},
\cite{Gonzalez-Kakehi-1}) 

When $\rank (G_{p,n})\,\leq\,\rank (G_{q,n})$, there is an
$O(n)$-invariant operator $\square_{p,q,n} : C^\infty
(G_{p,n})\,\rightarrow\,C^\infty (G_{p,n})$ which inverts $R_{p,q,n}$:
\begin{equation}\label{E:compactinversion}
f\;=\;\square_{p,q,n}\,R_{q,p,n}\,\circ R_{p,q,n}\,f,
\qquad\qquad\qquad\qquad f\in C^\infty (G_{p,n}).
\end{equation}

The operator $\square_{p,q,n}$, given explicitly in \cite{Kakehi},
corresponds to multiplication by a constant factor on each of the
$K$-types in $L^2 (G_{p,n})$, and is a differential operator when
$q-p$ is even.  We call $\square_{p,q,n}$ a {\it reproducing
operator}.

As stated in the introduction, we assume that $p<q$.  Let $\eta_0$ denote the $q$-plane
$\mathbb R e_1\,+\cdots\,+\,\mathbb R e_q$, and let $H_p$ and $H_q$
denote the subgroups of the Euclidean motion group
$M(n)\,=\,O(n)\rtimes \mathbb R^n$ fixing the $p$-plane $\sigma_0$ and
the $q$-plane $\eta_0$, respectively.  We have, in particular,
 $H_p\,=\,(O(p)\times
O(n-p))\rtimes \mathbb R^p$ and $G(p,n)\,=\,M(n)/H_p$.

For suitable compatible measures on $M(n),\;H_p,\;H_q$, and $H_p\cap
H_q$, the affine Grassmannian transform $\mathcal R^{(p,q)}$ is the Radon transform associated
with the double fibration
$$
\begin{matrix}
       &           &  M(n)/(H_p \cap H_q)   &          &     \\
       & \swarrow  &               & \searrow &     \\
  G(p,n)= M(n)/H_p  &           &               &          & G(q,n) = M(n)/H_q
\end{matrix}
$$
The corresponding incidence relation between $p$\,- and $q$\,-planes
in $\mathbb R^n$ is just inclusion.

We can write $\mathcal R^{(p,q)}$ explicitly in the following way
(\cite{Gonzalez-Kakehi-1}).  $G(p,n)$ is a vector bundle over
$G_{p,n}$ of rank $n-p$: its fiber over each $\sigma\in G_{p,n}$ is
$\sigma^\perp$; each $\ell\in G(p,n)$ is written uniquely as
$\ell\,=\,(\sigma,x)$, where $\sigma$ is the parallel translate of
$\ell$ through the origin and $\{x\}\,=\,\sigma^\perp\cap\ell$.  We
parametrize $G(q,n)$ in a similar manner: $G(q,n)\,=\,\{(\eta,v)\in
G_{q,n}\times \mathbb R^n\,|\,v\perp \eta\}$.  From
\cite{Gonzalez-Kakehi-1} the transform $\mathcal
R^{(p,q)}$ is then given by the formula
\begin{equation}
\mathcal R^{(p,q)}
f(\eta,v)\;=\;\int_{\sigma\subset\eta}\left(\int_{\sigma^\perp\cap\eta}
f(\sigma,v+x)\,dx
\right)\,d\sigma,\label{E:Rpq-expression}
\end{equation}
for any appropriate function $f$ on $G(p,n)$.  
The outer integral is taken over the set $\{\sigma\in
G_{p,n}\,|\,\sigma\subset\eta\}$, with respect to the normalized
measure invariant under all $u\in O(n)$ preserving $\eta$.

Let $\mathcal S(G(p,n))$ and $\mathcal S(G(q,n))$ denote the spaces of
Schwartz-class functions on $G(p.n)$ and $G(q,n)$, respectively.
Then $\mathcal R^{(p,q)}:\mathcal S(G(p,n))\rightarrow \mathcal
S(G(q,n))$, by \cite{Gonzalez-Kakehi-1}, \S 6.
 Next let $\mathcal F_p$ and $\mathcal F_q$ denote the {\it
partial Fourier transform} (i.e. Fourier transform on the fibers) on
$G(p,n)$ and $G(q,n)$, respectively:
\begin{equation}\qquad\qquad
\mathcal F_p f(\sigma,y)\;=\;\int_{\sigma^\perp}
f(\sigma,x)\,e^{-i\langle y,x\rangle}\,dx, \qquad\qquad\qquad\qquad
f\in\mathcal S(G(p,n)).
\end{equation}
Using (\ref{E:Rpq-expression}), it is not hard to prove the affine
Grassmannian version of the {\it
projection-slice theorem}:
\begin{equation}\label{E:projectionslice}
\mathcal F_q\circ \mathcal R^{(p,q)}
f(\eta,y)\;=\;\int_{\sigma\subset\eta} \mathcal F_p f(\sigma,y)\,d\sigma
\end{equation}
for all $(\eta,y)\in G(q,n)$.  (\cite{Gonzalez-Kakehi-1}, Prop. 6.1.)

We define the {\it rank} of $G(p,n)$ to be $\min(p+1,n-p)$.  
Since $\dim G(p,n) = (p+1)(n-p)$, we note that
$\rank(G(p,n))\,\leq\,\rank(G(q,n))$ if and only if
$\dim G(p,n) \,\leq\,\dim G(q,n) $.  The range
result
Theorem~\ref{T:diffrangechar} from \cite{Gonzalez-Kakehi-1} is essentially proved using the
projection-slice theorem and the range characterization of the compact
Radon transform $R_{p,q}$ on Grassmannians in $\mathbb R^{n-1}$ via
$O(n-1)$-invariant differential operators.  In addition, the
projection-slice theorem and the inversion formula
(\ref{E:compactinversion}) for $R_{p,q}$ on $\mathbb R^{n-1}$ are used to
prove an inversion formula  (Theorem 6.4 in \cite{Gonzalez-Kakehi-1}) for $\mathcal
R^{(p,q)}$ when $\dim G(p,n) \,=\,\dim G(q,n)$ and $q-p$ is even.
(See Rubin \cite{Rubin2} for another inversion formula for $\mathcal R^{(p,q)}$
without the parity restriction.)

It is not hard to obtain an analogue of the moment conditions (H) for the
transform $\mathcal R^{(p,q)}$.  Let $k\in\mathbb Z^+$ and
$f\in\mathcal S(G(p,n))$.  Then for any $(\eta,y)\in G(q,n)$ we have
\begin{equation}\label{E:momentderivation}
\begin{aligned}
\int_{\eta^\perp}\mathcal R^{(p,q)} f(\eta,v)\,{\langle
v,y\rangle}^k\,dv\;&=\;\int_{\eta^\perp}{\langle
v,y\rangle}^k\,\int_{\sigma\subset\eta}\,\int_{\sigma^\perp\cap\eta}
f(\sigma,v+x)\,dx\,d\sigma\,dv\\
&=\;\int_{\sigma\subset\eta}\,\int_{\eta^\perp}\,\int_{\sigma^\perp\cap\eta}
f(\sigma,v+x)\,{\langle v+x,y\rangle}^k\,dx\,dv\,d\sigma\\
&=\;\;\int_{\sigma\subset\eta}\,\int_{\sigma^\perp}
f(\sigma,w)\,{\langle w,y\rangle}^k\,dw\,d\sigma.
\end{aligned}
\end{equation}
Here we have used the fact that $\sigma^\perp$ equals the orthogonal
direct sum $(\sigma^\perp\cap\eta)\oplus\eta^\perp$.  The inner
integral above represents a smooth function $P_k(\sigma,y)$ on
$G(p,n)$:
\begin{equation}\label{E:momentpoly1}
\qquad\qquad P_k(\sigma,y)\;=\;\int_{\sigma^\perp}
f(\sigma,w)\,{\langle w,y\rangle}^k\,dw,\qquad\qquad\qquad\qquad y\in\sigma^\perp.
\end{equation}
Clearly $P_k$ is a homogeneous degree $k$ polynomial on the fibers of
$G(p,n)$. 

Equations (\ref{E:momentderivation}) and (\ref{E:momentpoly1}) lead us to define 
$\mathcal S_H(G(q,n))$ as the vector space consisting of all $\varphi\in
\mathcal S(G(q,n))$ satisfying the following condition 

\noindent\hspace{20pt}  {\it $(H')$: For each $k\in\mathbb Z^+$, there exists
a $C^\infty$ function $P_k$ on $G(p,n)$ such that}
\begin{enumerate}
\item {\it For any $\sigma\in G_{p,n}$, the function $y\mapsto
P_k(\sigma,y)$ is a homogeneous polynomial of degree $k$ on
$\sigma^\perp$.}
\item {\it For all $(\eta,y)\in G(q,n)$, we have}
\begin{equation}\label{E:P-integral}
\int_{\eta^\perp}\varphi(\eta,v)\,{\langle
v,y\rangle}^k\,dv\;=\;\int_{\sigma\subset\eta} P_k(\sigma,y)\,d\sigma.
\end{equation}
\end{enumerate}
\medskip
(\ref{E:momentderivation}) shows that the range $\mathcal
R^{(p,q)}\mathcal S(G(p,n))$ is a subspace of $\mathcal S_H(G(q,n))$.
Note that the condition ($H'$) reduces to the classical condition ($H$)
when $p=0$ and $q=n-1$.  We now state our main result:

\begin{theorem}\label{T:momentconditions}
Suppose that $p<q$ and $\rank(G(p,n))=\rank(G(q,n))$.  Then  $\mathcal
R^{(p,q)}\mathcal S(G(p,n))\,=\,\mathcal S_H(G(q,n))$.
\end{theorem}

Our proof roughly follows the lines of the classical proof. Let
$\varphi\in \mathcal S_H(G(q,n))$.  The partial Fourier transform
$\widetilde\varphi$ of $\varphi$,
\begin{equation}\label{E:definepft}
\widetilde\varphi(\eta,y)\;=\;\mathcal F_q\varphi(\eta,y)\;=
\;\int_{\eta^\perp}\varphi(\eta,v)\,e^{-i\langle
v,y\rangle}\,dv,
\end{equation}
belongs to $\mathcal S(G(q,n))$ by \cite{Gonzalez-Kakehi-1}.  We now
introduce the ``flag'' manifold $F_{q,n}=\{(\eta,\omega)\in
G_{q,n}\times S^{n-1}\,|\,\eta\perp\omega\}$.  (Define $F_{p,n}$
similarly.)  Then
$\widetilde\varphi$ gives rise to a smooth function $\widetilde\Phi$
on $F_{q,n}\times\mathbb R$
\begin{equation}\label{E:definePhi}
\widetilde\Phi(\eta,\omega,r)\;=\;\widetilde\varphi(\eta,r\omega).
\end{equation}
Note that
$\widetilde\Phi(\eta,\omega,r)=\widetilde\Phi(\eta,-\omega,-r)$.

For each $\omega\in S^{n-1}$, let $G_p(\omega^\perp)$ and
$G_q(\omega^\perp)$  denote the compact Grassmann manifolds of $p$-
and $q$-dimensional subspaces of the ($n-1$)-dimensional space
$\omega^\perp\subset\mathbb R^{n}$.  Then $G_p(\omega^\perp)$ and
$G_q(\omega^\perp)$ are diffeomorphic to $G_{p,n-1}$ and $G_{q,n-1}$,
respectively, and are homogeneous spaces of the subgroup $O(\omega)$
of $O(n)$ fixing $\omega$.

$F_{q,n}$ is a fiber bundle over $S^{n-1}$ with fibers
$G_p(\omega^\perp)$.  If we identify the Grassmannian $G_{q,n-1}$ with
$G_q(e_n^\perp)$, we can see that $F_{q,n}$ is the associated fiber
bundle $O(n)\times_{O(n-1)} G_{q,n-1}$ of the principal bundle
$O(n)\rightarrow S^{n-1},\quad u\mapsto u\cdot e_n$.  Let
$\widetilde\pi_q: O(n)\times G_{q,n-1}\rightarrow F_{q,n},\;
(u,\eta)\mapsto (u\cdot\eta, u\cdot e_n)$ be the quotient map.  Using
local cross sections, it is easy to see that a function $\Phi$ is
smooth on $F_{q,n}$ iff its lift $\Phi\circ\widetilde\pi_q$ is smooth
on $O(n)\times G_{q,n-1}$.

For each $\omega\in S^{n-1}$, let $R_{p,q}^\omega:
C^\infty(G_p(\omega^\perp))\rightarrow C^\infty(G_q(\omega^\perp))$ be
the $O(\omega)$-invariant Radon transform corresponding to the
inclusion incidence relation, and let $R_{q,p}^\omega$ be the dual
transform.  $R_{p,q}^\omega$ is of course just a translate, under
$O(n)$, of the transform $R_{p,q,n-1}:C^\infty(G_{p,n-1})\rightarrow
C^\infty(G_{q,n-1})$ defined in the beginning of this section.  Since
$p+q+1=n$, it follows that
$\rank(G_p(\omega^\perp))=\rank(G_q(\omega^\perp))=\min(p,q)$, and so 
 $R_{p,q}^\omega$ is a linear bijection.  Let
$\square_{p,q}^\omega$ be the corresponding reproducing operator; the
appropriate translate of (\ref{E:compactinversion}) for $R_{p,q}^\omega$
is
\begin{equation}\label{E:Rpqomegainversion}
\psi\;=\;\square_{p,q}^\omega R_{q,p}^\omega\circ R_{p,q}^\omega\,\psi,
\end{equation}
for all $\psi\in C^\infty(G_p(\omega^\perp))$.

Let us now return to the function $\widetilde\Phi$ on
$F_{q,n}\times\mathbb R$.  Since $R_{p,q}^\omega$ is a bijection,
there is, for each $(\omega,r)$,  a unique smooth function $\widetilde
F(\cdot,\omega;r)$ on $G_p^\omega$
 such that
\begin{equation}\label{E:flagint1}
\widetilde\Phi(\eta,\omega;r)\;=\;\int_{\sigma\subset\eta}
\widetilde F(\sigma,\omega;r)\,d\sigma,
\end{equation}
for all $(\eta,\omega;r)\in F_{p,n}\times\mathbb R$.  We also
express this as
$\widetilde\Phi(\eta,\omega;r)\,=\,\left(R_{p,q}^\omega \widetilde 
F (\cdot,\omega;r)\right)(\eta)$.  We can, of course, think of $\widetilde F$ as
being a function on $F_{p,n}\times\mathbb R$.

We want to prove that $\widetilde F$ is smooth on $F_{p,n}\times
\mathbb R$. (This is not completely obvious.)  Since the variable $r$
is fixed in (\ref{E:flagint1}),
this reduces to showing that the map $(\sigma,\omega)\mapsto\widetilde
F(\sigma,\omega;r)$ is $C^\infty$ on $F_{p,n}$ for each $r$.  In view
of the inversion formula
(\ref{E:Rpqomegainversion}), let us consider the integral transform $S$, from
functions on $F_{q,n}$ to functions on $F_{p,n}$, given by
\begin{equation}\label{E:flagint2}
S \Phi(\sigma,\omega)\;=\;
\int_{\sigma \subset \eta \subset \omega^{\perp}}
\Phi(\eta,\omega)\,d\eta\;=\;\left(R_{q,p}^\omega\Phi(\cdot,\omega)\right)(\sigma),
\qquad\qquad
(\sigma,\omega)\in F_{p,n}
\end{equation}
for all $\Phi\in C^\infty(F_{q,n})$.
Here $d \eta$ denotes the canonical and normalized measure
on the submanifold
$\{ \eta \in G_{q,n} \, | \, \sigma \subset \eta \subset \omega^{\perp} \}$
of $G_{q,n}$.

\begin{lemma}\label{T:S-smoothness}
$S$ is a continuous linear operator
from $C^\infty(F_{q,n})$ to $C^\infty(F_{p,n})$.
\end{lemma}
\begin{proof}
In fact $S$ is the Radon transform
associated with the double fibration
$$
\begin{matrix}
       &           &  O(n)/(\widetilde K_q \cap\widetilde K_p)   &          &     \\
       & \swarrow  &               & \searrow &     \\
  F_{q,n}= O(n)/\widetilde K_q  &           &               &
& F_{p,n} = O(n)/\widetilde K_p
\end{matrix}
$$
where $\widetilde K_q\,=\,O(q)\times O(n-q-1)$ and $\widetilde
K_p\,=\,O(p)\times O(n-p-1)$ are the subgroups of $O(n)$ fixing
$(\eta_0,e_n)\in F_{q,n}$ and $(\sigma_0,e_n)\in F_{p,n}$,
respectively.
Hence by Ch. I, Proposition 3.8 of \cite{GGA},
the transform $S$ is a continuous linear operator from $C^\infty(F_{q,n})$ to
$C^\infty(F_{p,n})$. 
\end{proof}

It will be useful to express $S$ in terms of associated fiber
bundles.  In terms of the quotient maps $\widetilde\pi_q: O(n)\times
G_{q,n-1}\rightarrow F_{q,n}$ and $\widetilde\pi_p: O(n)\times
G_{p,n-1}\rightarrow F_{p,n}$, we have
\begin{equation}\label{E:assocbundle1}
S\Phi\circ\widetilde\pi_p
(u,\sigma)\;=\;\int\limits_{\substack{\eta\supset\sigma\\ \eta\perp e_n}}\Phi\circ\widetilde\pi_q(u,\eta)\,d\eta\;=\;\left(R_{q,p,n-1}(\Phi\circ\widetilde\pi_q)(u,\cdot)\right)(\sigma).
\end{equation}
If $\Phi$ is smooth on $F_{q,n}$, then the right hand side is smooth
on $O(n)\times G_{p,n-1}$.  
(This also shows that $S\Phi$ is smooth on
$F_{p,n}$ whenever $\Phi$ is smooth in $F_{q,n}$.)

Next we define the operator $\square^{(p)}$ on $F_{p,n}$ by putting 
\begin{equation}\label{E:square-def}
\square^{(p)} F(\sigma,\omega)\,=\,\square_{p,q}^\omega
F(\cdot,\omega)(\sigma)
\end{equation}
for all $F\in C^\infty(F_{p,n})$.  (In the above, the operator
$\square_{p,q}^\omega$ 
acts on the first argument.)

\begin{lemma}\label{T:associatedbundle}
$\square^{(p)}$ is a continuous linear operator on
$C^{\infty} (F_{p,n})$.
\end{lemma}

\begin{proof}
Again letting $\widetilde\pi_p:O(n)\times G_{p,n-1}\rightarrow
F_{p,n}$ be the quotient map, the $O(n-1)$-invariance of
$\square_{p,q,n-1}$ implies that
\begin{equation}\label{E:square-action}
(\square^{(p)} G)\circ\widetilde\pi_p (u,\sigma)\;=\;\square_{p,q,n-1}\left(
G\circ\widetilde\pi_p\right)(u,\sigma),
\qquad \text{for} \; G \in C^\infty(F_{p,n})
\end{equation}
where $\square_{p,q,n-1}$ acts on the second argument.
Therefore, the continuity of $\square^{(p)}$ follows from (\ref{E:square-action})
and the continuity of $\square_{p,q,n-1}$.
\end{proof}

In particular, $\square^{(p)} F$ is a smooth function on $F_{p,n}$.

We now slightly modify the definitions of the operators $S$ and $\square^{(p)}$ in
(\ref{E:flagint2}) and (\ref{E:square-def}) so that they act on functions
on $F_{q,n}\times\mathbb R$ and $F_{p,n}\times\mathbb R$,
respectively.  In other words we put 
\begin{equation}\label{E:modify-S}
S \Psi(\sigma,\omega;r)\;=\;
\int_{\sigma \subset \eta \subset \omega^{\perp}} \, 
\Psi(\eta,\omega;r)\,d\eta,
\qquad\qquad\qquad\qquad
(\sigma,\omega)\in F_{p,n}
\end{equation}
for all $\Psi\in C^\infty(F_{q,n}\times\mathbb R)$, and
\begin{equation}\label{E:modify-square}
\square^{(p)} V(\sigma,\omega;r)\,=\,\square_{p,q}^\omega
V(\cdot,\omega;r)(\sigma)
\end{equation}
for all $V\in C^\infty(F_{p,n}\times\mathbb R)$.
Lemmas~\ref{T:S-smoothness} and \ref{T:associatedbundle}, suitably modified, still apply
to show that $S: C^\infty(F_{q,n}\times\mathbb R)\rightarrow
C^\infty(F_{p,n}\times\mathbb R)$ and 
$\square^{(p)}:C^\infty(F_{p,n}\times\mathbb R) \rightarrow
C^\infty(F_{p,n}\times\mathbb R)$
are continuous linear operators.

Now by the inversion formula
(\ref{E:Rpqomegainversion}) for $R_{p,q}^\omega$ and definitions (\ref{E:flagint2}) and 
(\ref{E:square-def}), we can recover $\widetilde F$ from
$\widetilde \Phi$ in equation~(\ref{E:flagint1}):
\begin{equation}\label{E:widetildeF}
\widetilde F(\sigma,\omega;r)\;=\;\left(\square^{(p)} \circ 
S\widetilde\Phi\right)(\sigma,\omega;r).
\end{equation}
By the remarks above, we see that $\widetilde F\in
C^\infty(F_{p,n}\times \mathbb R)$.

The uniqueness of $\widetilde F$ in (\ref{E:flagint1}) implies that 
$\widetilde F (\sigma,\omega;r)\,=\,\widetilde F (\sigma, -\omega;-r)$
for all $(\sigma,\omega;r)\in F_{p,n}\times \mathbb R$.  We next show
that the moment conditions ($H'$) for $k=0$ imply that 
$\widetilde F(\sigma,\omega;0)$ is constant in $\omega\in
S^{n-1}\cap\sigma^\perp$.
The function $P_0(\sigma,v)$ on $G(p,n)$ given in ($H'$) in this case is a
$0$th degree polynomial in $v$ for each $\sigma\in G_{p,n}$,
and thus depends only on $\sigma$,
so we put $P_0'(\sigma)\,=\,P_0(\sigma,0)$. We have $P_0'\in
C^\infty(G_{p,n})$ and by (\ref{E:P-integral}),
$$
\widetilde\phi(\eta,0)\;=\;\int_{\sigma\subset\eta}
P_0'(\sigma)\,d\sigma.
$$
If we take any $\omega\in S^{n-1}$ and $\eta\in G_q^\omega$, 
we have by (\ref{E:definePhi}) and (\ref{E:flagint1}),
$$
\widetilde\phi(\eta,0)\;=\;\widetilde\Phi(\eta,\omega,0)\;=
\;\int_{\sigma\subset\eta}\widetilde F(\sigma,\omega,0)\,d\sigma.
$$
Since $R_{p,q}^\omega$ is injective we see that $\widetilde
F(\sigma,\omega,0)\,=\,P_0'(\sigma)$ for all $\sigma\perp\omega$.

In view of this and the fact that $\widetilde F$ is even in $(\omega, r)$,
there exists a function $\widetilde f$ on $G(p,n)$ given by
\begin{equation}\label{E:define-f}
\qquad\qquad \widetilde f(\sigma,r\omega)\;=\;\widetilde
F(\sigma,\omega;r), \qquad\qquad\qquad\qquad (\sigma,\omega;r)\in
F_{p,n}\times \mathbb R.
\end{equation}
Now the mapping $(\sigma,\omega;r)\mapsto (\sigma,r\omega)$ is a local
diffeomorphism from $F_{p,n}\times(\mathbb R\setminus \{0\})$ onto
$G(p,n)\setminus G_{p,n}$.  (This is best seen by viewing both
$F_{p,n}\times \mathbb R$ and $G(p,n)$ as bundles over $G_{p,n}$, or
as associated bundles of the principal bundle $O(n)\rightarrow
G_{p,n},\;u\mapsto u\cdot \sigma_0$.)

$\widetilde f$ is therefore smooth on $G(p,n)\setminus G_{p,n}$.  In
addition, it is continuous on $G(p,n)$, since the map
$(\sigma,\omega;r)\mapsto (\sigma,r\omega)$ is a quotient map of
$F_{p,n}\times \mathbb R$ onto $G(p,n)$.  From
equation~(\ref{E:flagint1}), $\widetilde f$ satisfies the relation
\begin{equation}\label{E:flagint3}
\widetilde \phi(\eta,y)\;=\;\int_{\sigma\subset\eta} \widetilde f(\sigma,y)\;d\sigma
\end{equation}
for all $(\eta,y)\in G(q,n)$.

Our next objective, given in Propositions \ref{T:f-smoothness} and
\ref{T:f-rapid-decrease} below, is to prove that $\widetilde f$ is smooth on all of
$G(p,n)$, and that in fact $\widetilde f\in \mathcal S(G(p,n)$.  
Assuming this, let $f$ be the inverse partial Fourier transform of
$\widetilde f$:  $\mathcal F_p f=\widetilde f$.  Then the
projection-slice theorem 
\ref{E:projectionslice},
in conjunction with equations~(\ref{E:flagint1}), (\ref{E:definepft})
(\ref{E:definePhi}), and (\ref{E:define-f})
show that 
$$
\mathcal F_q\circ\mathcal
R^{(p,q)}f(\eta,r\omega)\,=\,\int_{\sigma\subset\eta} \widetilde
f(\sigma,r\omega)\,d\sigma\,=\,\mathcal F_q\varphi(\eta,r\omega)
$$
for all $(\eta,\omega)\in F_{q,n}$ and all $r\in\mathbb R$.
By the injectivity of $\mathcal F_q$, we get $\mathcal R^{(p,q)}
f=\varphi$, which proves Theorem~\ref{T:momentconditions}.

\begin{proposition}\label{T:f-smoothness} 
$\widetilde f\in C^\infty(G(p,n))$.
\end{proposition} 

\begin{proposition}\label{T:f-rapid-decrease} 
$\widetilde f\in\mathcal S(G(p,n))$.
\end{proposition}

We will give the proofs of the above two propositions
in Section \ref{Smoothness}.

\end{section}


\begin{section}{Differential operators on Grassmann manifolds and flag manifolds}\label{Diff-op}

In this section, we will study the calculus of differential operators
on $F_{p,n}$ and on 
$G(p,n) \setminus (G_{p,n} \times \{ 0 \}) \cong F_{p,n} \times \mathbb{R}_+$.
(Here we identify $G_{p,n} ( \subset G(p,n))$ with $G_{p,n} \times \{ 0 \}$,
using the parametrization 
$G(p,n) \ni \ell = (\sigma, x), \; \sigma \in G_{p,n}, x \in \sigma^{\perp}$.)
In particular, we will give a kind of {\it polar coordinate decomposition}
of differential operators on $G(p,n) \setminus (G_{p,n} \times \{ 0 \})$
analogous to (\ref{E:polardiff}).
The results in this section will be applied to the proofs of Propositions
\ref{T:f-smoothness} and \ref{T:f-rapid-decrease}
in Section \ref{Smoothness}.

Let $M(n) = O(n) \rtimes {\mathbb R}^n$ be the Euclidean motion group,
and let $\frak{m} (n) = so(n) \rtimes {\mathbb R}^n$ be its Lie algebra.
$\m (n)$ has basis consisting of the elementary skew symmetric matrices
$X_{i j}= E_{i j}-E_{j i} \; (1\leq i < j \leq n)$ and 
$E_k \; (1 \leq k \leq n)$, where $E_1, \cdots , E_n$ denote the infinitesimal
translations in the directions of the unit vectors $\mathbf{e}_1, \cdots , \mathbf{e}_n$.
Let $\u (\m(n) )$ denote the universal enveloping algebra
of $\m (n)$.
We note that $\u ({\mathbb R}^n) \subset \u (\m (n))$
is a commutative subalgebra of $\u ( \m (n))$.
We call $P \in \u ({\mathbb R}^n)$ a {\it homogeneous element}
if $P$ is written as a homogeneous polynomial of $E_1, \cdots , E_n$.
By the Poincar\'{e}-Birkhoff-Witt theorem,
the following proposition is easily obtained.

\begin{proposition}\label{Diff-Op-lemma:0}
Any element $U \in \u (\m (n))$ is written as
a linear combination of terms of the form
$Q P$, where $Q \in \u (so(n))$ and $P \in \u ({\mathbb R}^n)$
is a homogeneous element.
\end{proposition}

For the sake of simplicity, throughout this section, we write
the action of $X \in \m (n)$ on $f \in C^{\infty} (G(p,n)$ as
\begin{equation*}
X \cdot f ((\sigma, x)) = \frac{d}{dt} f(\exp (-tX) \cdot (\sigma, x)) |_{t=0}.
\end{equation*}
Similarly, if $X \in so(n)$ and $f \in C^{\infty} (F_{p,n})$
(or $f \in C^{\infty} (G_{p,n})$ ), we write the action of $X$ on $f$
as $X \cdot  f$.
In addition, we also write the action of $U \in \u (\m (n))$
on $f \in C^{\infty} (G(p,n)$ as $U \cdot  f$.

For $v \in \mathbf{S}^{n-1}$ and $\alpha >0$,
we introduce the open sets
\begin{equation*}
F_{p,n}^{\alpha} (v) \; = \;
 \{ \; (\sigma, \omega) \in F_{p,n} \; | \;
   || \, \Pr_{\sigma^{\perp}} \, v \, || > \alpha \; \},
\end{equation*}
$\Pr_{\sigma^{\perp}}$ denoting orthogonal projection to $\sigma^{\perp}$.
In addition, for an index set 
$I = \{ \; i_1, \cdots , i_m \; \}$ with $1 \leq i_1 < \cdots < i_m \leq n$,
let
\begin{equation*}
F_{p,n}^{\alpha} [I] = \bigcap_{k=1}^{m} \, F_{p,n}^{\alpha} (\mathbf{e}_{i_k}).
\end{equation*}
Then we have the following.

\begin{lemma}\label{Diff-Op-lemma:1}
\begin{equation*}
\bigcup_{\# I = n-p} \, F_{p,n}^{\alpha} [I]  \; = \; F_{p,n},
\qquad \text{for some} \; \; \alpha >0.
\end{equation*}
\end{lemma}

\begin{proof} \quad
Take an arbitrary point $(\sigma, \omega) \in F_{p,n}$.
Since $\sigma$ is a $p$-dimensional subspace of $\mathbb{R}^n$,
there exist $n-p$ unit vectors $\mathbf{e}_{i_1}, \cdots, \mathbf{e}_{i_{n-p}}$
($1 \leq i_1 < \cdots < i_{n-p} \leq n$) such that 
$\mathbf{e}_{i_k} \notin \sigma, \; (1 \leq k \leq n-p)$.
So if we put
\begin{equation*}
\alpha = \frac{1}{2} \, \min 
  \{ \, || \Pr_{\sigma^{\perp}} \, \mathbf{e}_{i_k} || \; | \; 1 \leq k \leq n-p \; \},
\quad I = \{ \; i_1, \cdots , i_{n-p} \; \},
\end{equation*}
then $(\sigma, \omega) \in F_{p,n}^{\alpha} [I]$.
Let $V_{\alpha} = \bigcup_{\# I = n-p} \; F_{p,n}^{\alpha} [I]$.
Then, the above result shows that $\bigcup_{\alpha >0} \, V_{\alpha} \; = F_{p,n}$.
Since $\{ \, V_{\alpha} \, \}_{\alpha >0}$ is an open covering of the compact set
$F_{p,n}$, we can take a finite covering $\{ \, V_{\alpha_j} \, \}_{j=1}^N$
of $F_{p,n}$.
Let $\alpha_0 = \min_{1 \leq j \leq N} \, \alpha_j.$ Then we have
$V_{\alpha_0} =F_{p,n}$, which proves the assertion.
\end{proof}

\begin{lemma}\label{Diff-Op-lemma:2}
Let 
$\mathcal{U}_N = \{ \; \theta = (\theta_1, \cdots, \theta_N) 
  \in \mathbf{S}^{N-1} \; | \; \theta_N >0 \; \}$.
There exist smooth functions $a_j \; (j=1, \cdots, N-1)$ and $b$ on
$\mathbf{S}^{N-1} \times \mathcal{U}_N$ such that
\begin{equation*}
\begin{split}
( E_v \, f ) ( r \theta) \; &= \; \sum_{j=1}^{N-1} \;
  \frac{a_j (v, \theta)}{r} \, X_{j N} f( r \theta) 
    + b (v, \theta) \frac{\partial}{\partial r} f (r \theta),\\
\quad &\text{for} \; (v, \theta) \in \mathbf{S}^{N-1} \times \mathcal{U}_N, \; r>0,
\quad \text{and for} \; f \in C^{\infty} (\mathbb{R}^N \setminus \{ 0 \} ).
\end{split}
\end{equation*}
Here $E_v$ denotes the directional derivative in the direction of
$-v \in \mathbb{R}^N$, namely,
\begin{equation*}
E_v \cdot f (x) = \frac{d}{dt} f(x -tv) |_{t=0}, \qquad \text{for} \;
f \in C^{\infty} (\mathbb{R}^N).
\end{equation*}
\end{lemma}

\begin{proof}
It is not hard to see that
\begin{equation*}
\frac{1}{r} X_{j N} 
 = \theta_N \, \frac{\partial}{\partial x_j} 
   - \theta_j \, \frac{\partial}{\partial x_N},
\quad
\frac{\partial}{\partial r} = \sum_{j=1}^N \; \theta_j \, \frac{\partial}{\partial x_j} .
\end{equation*}
Hence
\begin{equation*}
\begin{pmatrix}
\frac{1}{r} X_{1 N}\\
\vdots             \\
\frac{1}{r} X_{N-1 N}\\
\frac{\partial}{\partial r}
\end{pmatrix}
\; = \;
\begin{pmatrix}
\theta_N    &     0    &    \cdots     &       0       & -\theta_1     \\
    0       & \theta_N &     0         &    \cdots     & -\theta_2     \\
 \vdots     &          &    \ddots     &               &  \vdots       \\
    0       & \cdots   &      0        &  \theta_N     & -\theta_{N-1} \\
\theta_1    & \theta_2 &    \cdots     &  \theta_{N-1} &  \theta_N            
\end{pmatrix}
\;
\begin{pmatrix}
\frac{\partial}{\partial x_1}     \\
\vdots  \\
\frac{\partial}{\partial x_{N-1}}  \\
\frac{\partial}{\partial x_N}     \\
\end{pmatrix}.
\end{equation*}
Note that
\begin{equation*}
\det \;
\begin{pmatrix}
\theta_N    &     0    &    \cdots     &       0       & -\theta_1     \\
    0       & \theta_N &     0         &    \cdots     & -\theta_2     \\
 \vdots     &          &    \ddots     &               &  \vdots       \\
    0       & \cdots   &      0        &  \theta_N     & -\theta_{N-1} \\
\theta_1    & \theta_2 &    \cdots     &  \theta_{N-1} &  \theta_N   
\end{pmatrix}
\; = \; \theta_N^{N-2} >0, \qquad 
\text{if} \; \theta \in \mathcal{U}_N.
\end{equation*}
Thus, each vector field $\frac{\partial}{\partial x_j}$
is written in the form,
\begin{equation*}
\frac{\partial}{\partial x_j} = \sum_{k=1}^{N-1} \; 
   \frac{a_{jk} (\theta)}{\theta_N^{N-2}} \,
     \frac{1}{r} X_{k N}
    + \frac{b_j (\theta)}{\theta_N^{N-2}} \, \frac{\partial}{\partial r}, 
\quad (j=1, \cdots, N),
\end{equation*}
where $a_{jk} (\theta)$ and $b_j (\theta)$ are polynomials of
$\theta \in \mathcal{U}_N$.
Since $E_v = - \sum_{j=1}^N \, v_j \frac{\partial}{\partial x_j}$, the assertion follows
from the above expression of $\frac{\partial}{\partial x_j}$.
\end{proof}

Next, in a similar way to the case of ${\mathbb R}^N$
we introduce the radial derivative $E_r$ on 
$G(p,n) \setminus (G_{p,n} \times \{ 0 \} )$ 
and the directional derivative $E_v$ on $G(p,n)$ in the direction of
$-v \in \mathbb{R}^n$ as follows.
\begin{align*}
E_r f (\sigma, x) &= \frac{1}{||x||} \frac{d}{dr} f ( \sigma, rx) |_{r=1}, 
  \qquad \text{for} \; f \in 
    C^{\infty} (G(p,n) \setminus (G_{p,n} \times \{ 0 \} )),\\
E_v \cdot f (\ell) &= \frac{d}{dt} f(\ell -tv) |_{t=0}, \qquad \text{for} \;
f \in C^{\infty} (G(p,n) \setminus (G_{p,n} \times \{ 0 \})).
\end{align*}
Our aim now is to generalize (\ref{E:polardiff}) to $G(p,n)$.

\begin{proposition}\label{Diff-Op-proposition:1}
For an arbitrary point $(\sigma_1, \omega_1) \in F_{p,n}$, there exist
an open neighborhood $\mathcal{U}$ of $(\sigma_1, \omega_1)$ in $F_{p,n}$
and smooth functions $a_{i j}^l (\sigma, \omega)$,
$b^l (\sigma, \omega)$ on $\mathcal{U}$ such that
\begin{equation}\label{Diff-Op-proposition:1-1}
E_l  \cdot f( \sigma, r \omega) = \sum_{1 \leq i < j \leq n} \;
   \frac{1}{r} a_{i j}^l (\sigma, \omega) \, X_{ij}  \cdot f ( \sigma, r \omega) 
     + b^l (\sigma, \omega) \, E_r f ( \sigma, r \omega),
\end{equation}
for $l \; (1 \leq l \leq n)$, $(\sigma, \omega) \in \mathcal{U}$,
$r >0$ and for 
$f \in C^{\infty} \, (G(p,n) \setminus G_{p,n} \times \{ 0 \} )$.
\end{proposition}

\begin{proof}
Let us first take an open neighborhood $\mathcal{V}_1$ of $\sigma_1$
in $G_{p,n}$ and a smooth local cross section
$u : \mathcal{V}_1 \to SO(n)$ such that
\begin{equation*}
u(\sigma) \sigma = 
  \mathbb{R} \mathbf{e}_1 \oplus \cdots \oplus \mathbb{R} \mathbf{e}_p
 \; \; (\text{ we put } \equiv \sigma_0 ), \quad
\text{for} \; \sigma \in \mathcal{V}_1, \quad
u(\sigma_1) \omega_1 = \mathbf{e}_n.
\end{equation*}
Then, by taking a sufficiently small open neighborhood $\Omega_1$ of
$\omega_1$ in $\mathbf{S}^{n-1}$, we have
\begin{equation*}
u(\sigma) \omega \in 
\{ \; ( 0, \cdots , 0, \theta_1, \cdots, \theta_{n-p} )
  \in \mathbf{S}^{n-1} \; | \; \theta_{n-p} >0 \; \},
\end{equation*}
for $(\sigma, \omega) \in ( \mathcal{V}_1 \times \Omega_1 ) \cap F_{p,n}$.
Now, let
\begin{equation*}
\mathcal{U} = ( \mathcal{V}_1 \times \Omega_1 ) \cap F_{p,n}.
\end{equation*}
By Lemma \ref{Diff-Op-lemma:1}, there exists an index set 
$I = \{ \; i_1, \cdots , i_{n-p} \; \}$ such that
$\mathcal{U} \subset F_{p,n}^{\alpha} [I]$.
(Replace $\mathcal{V}_1$ and $\Omega_1$ by smaller neighborhoods
if necessary.)

\underline{Step 1.} \quad
Let us first consider the case $l  \in I$.
\par
The vector $\mathbf{e}_l$ can be written as
\begin{equation}\label{Diff-Op-proposition:1-3}
\begin{split}
\mathbf{e}_l &= \Pr_{\sigma} \mathbf{e}_l 
               + R_l (\sigma) \, y_l (\sigma),\\
R_l (\sigma) &= || \Pr_{\sigma^{\perp}} \mathbf{e}_l ||, \quad
  y_l (\sigma) = 
    \frac{\Pr_{\sigma^{\perp}} \mathbf{e}_l}{|| \Pr_{\sigma^{\perp}} \mathbf{e}_l ||}
      \in \mathbf{S}^{n-1} \cap \sigma^{\perp}.
\end{split}
\end{equation}
By the definition of $F_{p,n}^{\alpha} [I]$,
\begin{equation*}
R_l (\sigma) \equiv || \Pr_{\sigma^{\perp}} \mathbf{e}_l || > \alpha,
\quad \text{for} \; \sigma \in \mathcal{V}_1.
\end{equation*}
Therefore, $R_l (\sigma)$ and $ y_l (\sigma)$ are smooth on $\mathcal{V}_1$.
Moreover, by (\ref{Diff-Op-proposition:1-3})
\begin{equation*}
E_l  \cdot  f(\sigma, r \omega) = R_l (\sigma) E_{y_l (\sigma) }  \cdot  f(\sigma, r \omega),
\quad \text{for} \; \; 
f \in C^{\infty} \, (G(p,n) \setminus G_{p,n} \times \{ 0 \} ).
\end{equation*}
From now on, we will decompose the above vector field $E_{y_l (\sigma) }$ as
a linear combination of rotational derivatives and the radial derivative $E_r$.
Since 
$u(\sigma) \cdot y_l (\sigma) \in \sigma_0^{\perp}
 \, = \mathbb{R} \mathbf{e}_{p+1} \oplus \cdots \oplus \mathbb{R} \mathbf{e}_n
\cong \mathbb{R}^{n-p}$, we can apply Lemma \ref{Diff-Op-lemma:2}
to the case when $N= n-p$ and $v= u(\sigma) \cdot y_l (\sigma)$.
Thus
\begin{equation}\label{Diff-Op-proposition:1-4}
\begin{split}
E_{u(\sigma) \cdot y_l (\sigma)}  \cdot  g( r \theta) 
\;  = & \; 
\sum_{m=p+1}^{n-1} \; 
  \frac{1}{r} \, A_m (u(\sigma) \cdot y_l (\sigma), \theta) \,
    ( X_{m n}  \cdot  g) (r \theta)\\
&+ \; B (u(\sigma) \cdot y_l (\sigma), \theta) \,
   ( E_r \, g) (r \theta),
\end{split}
\end{equation}
for $\theta = (0, \cdots, 0, \theta_1, \cdots, \theta_{n-p} )
   \in \sigma_0^{\perp} \cap \mathbf{S}^{n-1}$, $\theta_{n-p} >0, \; \; r>0$,
and for $g \in C^{\infty} (\mathbb{R}^{n-p} \setminus \{ 0 \} )$.
Here $A_m$ and $B$ are smooth functions on
$\mathbf{S}^{n-p-1} \times 
\{ \; \theta = (0, \cdots, 0, \theta_1, \cdots, \theta_{n-p} )
\in \sigma_0^{\perp} \cap \mathbf{S}^{n-1} \; | \; \theta_{n-p} >0 \; \}$.
(We identify $\theta$ with a unit vector $(\theta_1, \cdots, \theta_{n-p} )$
in $\mathbb{R}^{n-p}$.)

Let us take
$g(r \theta) = f (u(\sigma)^{-1} \cdot \sigma_0, \, r u(\sigma)^{-1} \cdot \theta)$
in (\ref{Diff-Op-proposition:1-4}).
Then,
\begin{equation}\label{Diff-Op-proposition:1-5}
\begin{split}
E_{y_l (\sigma) }  \cdot  f(\sigma, r \omega)
 &= E_{u(\sigma) \cdot y_l (\sigma)}  \cdot 
    f (u(\sigma)^{-1} \cdot \sigma_0, \, r u(\sigma)^{-1} \cdot \theta)
      |_{\theta = u(\sigma) \cdot \omega}\\
 &= \sum_{m=p+1}^{n-1} \; 
   \frac{1}{r} \, A_m (u(\sigma) \cdot y_l (\sigma), \theta)  \,
       X_{m n}  \cdot f (u(\sigma)^{-1} \cdot \sigma_0, \, r u(\sigma)^{-1} \cdot \theta)
         |_{\theta = u(\sigma) \cdot \omega}\\
 & \; \; + B (u(\sigma) \cdot y_l (\sigma), \theta) \,
     E_r \, f (u(\sigma)^{-1} \cdot \sigma_0, \, r u(\sigma)^{-1} \cdot \theta)
         |_{\theta = u(\sigma) \cdot \omega}\\
 &= \sum_{m=p+1}^{n-1} \; 
   \frac{1}{r} \, A_m (u(\sigma) \cdot y_l (\sigma), \theta) \,
     (X_{m n}^{\tau (u(\sigma))}  \cdot  f) (\sigma, r \omega)\\
 & \; \; + B (u(\sigma) \cdot y_l (\sigma), \theta) \,
    (E_r^{\tau (u(\sigma))} \, f) (\sigma, r \omega).
\end{split}
\end{equation}
We note that the radial derivative $E_r$ is invariant under the action of $SO(n)$,
in particular, $E_r^{\tau (u(\sigma))} =E_r$.
In addition, we also note that $X_{m n}^{\tau (u(\sigma))}$ is written in the form
\begin{equation}\label{Diff-Op-proposition:1-6}
X_{m n}^{\tau (u(\sigma))} 
  = \sum_{1 \leq i < j \leq n} \; C_{i j}^m (\sigma) \, X_{i j},
\end{equation}
where $C_{i j}^m (\sigma)$ is a smooth function on $\mathcal{V}_1$.
Combining (\ref{Diff-Op-proposition:1-5}) and (\ref{Diff-Op-proposition:1-6}),
we have the following expression
\begin{equation}\label{Diff-Op-proposition:1-7}
E_{y_l (\sigma) }  \cdot  f(\sigma, r \omega)
= \sum_{1 \leq i < j \leq n} \;
   \frac{1}{r} \widetilde{\, a \,}_{i j}^l (\sigma, \omega) \, 
   X_{ij}  \cdot f ( \sigma, r \omega) 
     + \widetilde{\, b \,}^l (\sigma, \omega) \, E_r f ( \sigma, r \omega),
\end{equation}
for some smooth functions
$\widetilde{\, a \,}_{i j}^l$ and $\widetilde{\, b \,}^l$ on $\mathcal{U}$.
Since $R_l (\sigma)$ is smooth on $\mathcal{V}_1$, we obtain an expression
of the form (\ref{Diff-Op-proposition:1-1}) for 
$l \in I$.

\underline{Step 2.} \quad
Next, let us consider the case $l  \notin I$.
Since $\mathcal{U} \subset F_{p,n}^{\alpha} [I]$,
we see easily that
\begin{equation}\label{Diff-Op-proposition:1-8}
\mathbb{R}^n = 
  \sigma \oplus \langle \mathbf{e}_{i_1}, \cdots, \mathbf{e}_{i_{n-p}} \rangle,
    \quad \text{for} \; \; \sigma \in \mathcal{V}_1,
\end{equation}
where $\langle \mathbf{e}_{i_1}, \cdots, \mathbf{e}_{i_{n-p}} \rangle$ denotes
the $(n-p)$-dimensional subspace spanned by $\mathbf{e}_{i_1}, \cdots, \mathbf{e}_{i_{n-p}}$.
(Note that the decomposition (\ref{Diff-Op-proposition:1-8}) is not necessarily
an orthogonal decomposition.)
Using the decomposition (\ref{Diff-Op-proposition:1-8}),
we can write $\mathbf{e}_l$ in the form
\begin{equation*}
\mathbf{e}_l = z + \sum_{k=1}^{n-p} \, a_k (\sigma) \mathbf{e}_{i_k},
 \quad z \in \sigma.
\end{equation*}
Obviously the coefficient $a_k (\sigma)$ is a smooth function of $\sigma \in \mathcal{V}_1$.
Thus,
\begin{equation}\label{Diff-Op-proposition:1-10}
E_l \cdot f (\sigma, r \omega) =
  \sum_{k=1}^{n-p} \, a_k (\sigma) E_{i_k} \cdot f (\sigma, r \omega).
\end{equation}
From (\ref{Diff-Op-proposition:1-10}) and the result of Step 1, we can conclude that
$E_l f (\sigma, r \omega) = E_{\mathbf{e}_l} f (\sigma, r \omega)$ has an expression
of the form (\ref{Diff-Op-proposition:1-1}) for $l \notin I$.
\end{proof}

Proposition \ref{Diff-Op-proposition:1} yields the following.

\begin{proposition}\label{Diff-Op-proposition:2} \quad
There exist smooth functions $a_{i j}^l (\sigma, \omega)$ and
$b^l (\sigma, \omega)$ on $F_{p,n}$ such that
\begin{equation}\label{Diff-Op-proposition:2-1}
E_l \cdot f( \sigma, r \omega) = \sum_{1 \leq i < j \leq n} \;
   \frac{1}{r} a_{i j}^l (\sigma, \omega) \, X_{ij} \cdot f ( \sigma, r \omega) 
     + b^l (\sigma, \omega) \, E_r f ( \sigma, r \omega),
\end{equation}
for $l \; (1 \leq l \leq n)$, $(\sigma, \omega) \in F_{p,n}$,
$r >0$ and for 
$f \in C^{\infty} \, (G(p,n) \setminus G_{p,n} \times \{ 0 \} )$.
\end{proposition}

\begin{proof}
Since $F_{p,n}$ is compact, it follows from Proposition \ref{Diff-Op-proposition:1}
that there exists a finite open covering $\{ \mathcal{U}_{\nu} \}_{1 \leq \nu \leq m}$
of $F_{p,n}$ such that $E_l f$ has an expression
of the form (\ref{Diff-Op-proposition:1-1}) on each $\mathcal{U}_{\nu}$.
So the above (globally defined) smooth functions $a_{i j}^l (\sigma, \omega)$ and
$b^l (\sigma, \omega)$ can be constructed, using a partition of unity
for the open covering $\{ \mathcal{U}_{\nu} \}_{1 \leq \nu \leq m}$.
\end{proof}

Let us consider $\frak{u} (\mathbb{R}^n)$ to be
a subalgebra of 
$\frak{u} (\frak{m}(n)) = \frak{u} (\frak{so}(n) \oplus \mathbb{R}^n) $.
We see easily that an element of $\frak{u} (\mathbb{R}^n)$
is written as a polynomial of vector fields $E_1, \cdots, E_n$.
The following theorem is the affine Grassmann version of equation (\ref{E:polardiff}).

\begin{theorem}\label{Diff-Op-theorem:1}
Let $P = P^{(m)} (E_1, \cdots, E_n) \in \frak{u} (\mathbb{R}^n)$ be 
a homogeneous element of order $m$.
Then for 
$f \in C^{\infty} ( G(p, n) \setminus G_{p,n} \times \{ 0 \} )$,
$P^{(m)} (E_1, \cdots, E_n) \cdot  f ( \sigma, r \omega)$ is expressed as
\begin{equation*}
P^{(m)} (E_1, \cdots, E_n) \cdot  f ( \sigma, r \omega) \; = \;
 \sum_{k=0}^m \; A_k^{(m)} (D_{(\sigma, \omega)}) r^{-k} E_r^{m-k} \, f( \sigma, r \omega),
 \quad \text{for} \; \; (\sigma, \omega) \in F_{p,n}, \;
r>0,
\end{equation*}
where $A_k^{(m)} (D_{(\sigma, \omega)})$ is a differential operator
on $F_{p,n}$ of order at most $k$.
\end{theorem}

\begin{proof}
We will prove the theorem by induction with respect to $m$.
In the case $m=1$, the assertion follows from Proposition \ref{Diff-Op-proposition:2}.
Suppose that the assertion of the theorem holds for any element of 
$\frak{u} (\mathbb{R}^n)$ of order at most $m$.
Let us take any homogeneous
element $Q \in \frak{u} (\mathbb{R}^n)$ of order $m+1$.
Then, $Q$ can be written as a linear combination of
$E_j Q^{(j)} \; (1 \leq j \leq n)$, where $Q^{(j)} \in \frak{u} (\mathbb{R}^n)$
is a homogeneous element of order $m$.
Therefore, without loss of generality we may assume that $Q = E_1 P$
for some homogeneous element $P \in \frak{u} (\mathbb{R}^n)$ of order $m$.
Then, by the hypothesis of induction,
$P f ( \sigma, r \omega)$ is written as
\begin{equation*}
P (E_1, \cdots, E_n) \, f ( \sigma, r \omega) \; = \;
 \sum_{k=0}^m \; A_k (D_{(\sigma, \omega)}) r^{-k} E_r^{m-k} \, f( \sigma, r \omega)
 \quad \text{for} \; \; (\sigma, \omega) \in F_{p,n}, \;
r>0,
\end{equation*}
where $A_k (D_{(\sigma, \omega)})$ is a differential operator
on $F_{p,n}$ of order at most $k$.

On the other hand, by Proposition \ref{Diff-Op-proposition:2},
there exist smooth functions 
$a_{i j} (\sigma, \omega)$ and
$b (\sigma, \omega)$ on $F_{p,n}$ such that
\begin{equation*}
E_1 \cdot F( \sigma, r \omega) = \sum_{1 \leq i < j \leq n} \;
   \frac{1}{r} a_{i j} (\sigma, \omega) \, X_{ij} \cdot F ( \sigma, r \omega) 
     + b (\sigma, \omega) \, E_r F ( \sigma, r \omega),
\end{equation*}
for $(\sigma, \omega) \in F_{p,n}$,
$r >0$ and for 
$F \in C^{\infty} \, (G(p,n) \setminus G_{p,n} \times \{ 0 \} )$.
Here we note that
\begin{equation*}
r X_{ij} = X_{ij} r, \; \;
E_r X_{ij} = X_{ij} E_r, \; \;
r A_k (D_{(\sigma, \omega)}) = A_k (D_{(\sigma, \omega)}) r, \; \; 
E_r A_k (D_{(\sigma, \omega)}) = A_k (D_{(\sigma, \omega)}) E_r.
\end{equation*}
Thus we have
\begin{equation*}
\begin{split}
Q \cdot f ( \sigma, r \omega) 
&= E_1 P \cdot  f( \sigma, r \omega)\\
&= \sum_{1 \leq i < j \leq n} \; \sum_{k=0}^m \;
  a_{i j} (\sigma, \omega) X_{ij} \, A_k (D_{(\sigma, \omega)})
   r^{-k-1} E_r^{m-k} \, f( \sigma, r \omega)\\
& \; \; + \sum_{k=0}^m \;
  b (\sigma, \omega) \, A_k (D_{(\sigma, \omega)}) \,
    r^{-k} E_r^{m+1-k} \, f( \sigma, r \omega)\\
&= b (\sigma, \omega) \, A_0 (D_{(\sigma, \omega)}) \,
    E_r^{m+1} \, f( \sigma, r \omega)\\
& \; \; + \sum_{k=1}^m \;
     \left\{ \; b (\sigma, \omega) \, A_k (D_{(\sigma, \omega)})
       + \sum_{1 \leq i < j \leq n} \, 
         a_{i j} (\sigma, \omega) X_{ij} \, A_{k-1} (D_{(\sigma, \omega)})
           \right\} \; r^{-k} E_r^{m+1-k} \, f( \sigma, r \omega)\\
& \; \; + \sum_{1 \leq i < j \leq n} \; 
   a_{i j} (\sigma, \omega) X_{ij} \, A_m (D_{(\sigma, \omega)})
    \, r^{-m-1} \,  f( \sigma, r \omega).
\end{split}
\end{equation*}
It follows from the above expression that
the assertion of the theorem  holds for $Q = E_1 P$.
Therefore, the proof is completed.
\end{proof}

\end{section}


\begin{section}{Smoothness of ${\widetilde f}$}\label{Smoothness}

In this section, we will prove Propositions 
\ref{T:f-smoothness} and
\ref{T:f-rapid-decrease}.

\begin{proposition}\label{smoothness-proposition:1}
There exist homogeneous polynomials $P_k (\sigma, x)$ of degree $k$
on $G(p,n)$ ($k = 0, 1, 2, \cdots,$) such that
\begin{align*}
{\widetilde f}(\sigma, x) 
 &= \sum_{k=0}^{N-1} \; \frac{(-i)^k}{k!} \, P_k (\sigma, x)
  + (\square^{(p)} \, \Phi_N ) (\sigma, \omega ; r), \\
\text{where} \quad
\Phi_N (\sigma, \omega ; r)
 &= \int_{\sigma \subset \eta \subset \omega^{\perp}} \;
   \int_{y \in \eta^{\perp}} \;
     e_N (-i \langle y, r \omega \rangle ) \varphi (\eta, y)
        \, dy d \eta,
\end{align*}
for $N, \; (N  = 1, 2, 3, \cdots, )$
and for $(\sigma, x) \in G(p,n)$ with 
$x = r \omega, \; ((\sigma, \omega) \in F_{p,n})$.
Here, as in Section2, $e_N (t)$ denotes the $N$th remainder term of the Taylor expansion
of $e^t$.
Moreover, as in (\ref{E:flagint2}), $d \eta$ denotes the canonical and normalized measure
on the set 
$\{ \eta \in G_{q,n} \, | \, \sigma \subset \eta \subset \omega^{\perp} \}$.
\end{proposition}

\begin{proof}
By the definition of ${\widetilde F}$,
\begin{equation*}
\begin{split}
{\widetilde f} (\sigma, r \omega) 
&= \square_{p, q}^{\omega} \int_{\sigma \subset \eta \subset \omega^{\perp}} \;
      \mathcal{F}_q \varphi (\eta, r \omega) \, d \eta\\
&= \square_{p, q}^{\omega} \int_{\sigma \subset \eta \subset \omega^{\perp}}
     \int_{y \in \eta^{\perp}} \;
       e^{-i \langle y, r \omega \rangle } \,
         \varphi (\eta, y) \, dy d \eta\\
&= \square_{p, q}^{\omega} \int_{\sigma \subset \eta \subset \omega^{\perp}}
     \int_{y \in \eta^{\perp}} \;
       \left\{
         \sum_{k=0}^{N-1} \; \frac{(-i)^k}{k!} \, \langle y, r \omega \rangle^k
           + e_N (-i \langle y, r \omega \rangle )
             \right\}
           \varphi (\eta, y) \, dy d \eta\\
&= \sum_{k=0}^{N-1} \; \frac{(-i)^k}{k!} \, 
     \square_{p, q}^{\omega} \int_{\sigma \subset \eta \subset \omega^{\perp}}
       \int_{y \in \eta^{\perp}} \;
         \varphi (\eta, y)  \, \langle y, r \omega \rangle^k
           \, dy d \eta\\
& \; \; + (\square^{(p)} \, \Phi_N ) (\sigma, \omega ; r).
\end{split}
\end{equation*}
Since $\varphi$ satisfies the moment condition $(H')$,
there exist homogeneous polynomials $P_k (\sigma, x)$ of degree $k$
on $G(p,n)$ ($k = 0, 1, 2, \cdots $), such that
\begin{equation*}
\int_{\sigma \subset \eta} \;
 P_k (\sigma, r \omega) \; d \sigma
  \; = \;
   \int_{y \in \eta^{\perp}} \;
     \varphi (\eta, y) \, \langle y, r \omega \rangle^k \; dy,
\quad \text{for} \; \forall (\eta, r \omega) \in G(q, n) 
 \; \text{with} \; \eta \perp \omega.
\end{equation*}
Applying the inversion formula for $R_{p,q}^{\omega}$ to the both sides of the above
equality, we have
\begin{equation*}
P_k (\sigma, r \omega)
 = \square_{p, q}^{\omega} \int_{\sigma \subset \eta \subset \omega^{\perp}}
       \int_{y \in \eta^{\perp}} \;
         \varphi (\eta, y)  \, \langle y, r \omega \rangle^k
           \, dy d \eta.
\end{equation*}
Therefore, using these homogeneous polynomials $P_k, \; (k=0,1,2,\cdots,)$,
${\widetilde f} (\sigma, r \omega)$ is written as
\begin{equation*}
{\widetilde f} (\sigma, r \omega) 
 = \sum_{k=0}^{N-1} \; \frac{(-i)^k}{k!} \, P_k (\sigma, r \omega)
  + (\square^{(p)} \, \Phi_N ) (\sigma, \omega ; r),
\end{equation*}
which completes the proof.
\end{proof}

Since $G(p,n) \setminus (G_{p,n} \times \{ 0 \}) \cong F_{p,n} \times \mathbb{R}_+$,
the radial derivative $E_r$ also acts on a function on
$F_{p,n} \times \mathbb{R}_+$.
From now on, we extend the action of $E_r$ to 
$F_{p,n} \times (\mathbb{R} \setminus \{ 0 \})$ so that
$$
E_r \Psi (\sigma, \omega; r) = E_r \Psi (\sigma, -\omega; -r)
\quad \text{for} \; \; 
\Psi \in C^{\infty} (F_{p,n} \times (\mathbb{R} \setminus \{ 0 \})).
$$

\begin{proposition}\label{smoothness-proposition:2}
For $k, l \in \mathbb{Z}^+$ with $k+l=m$,
we have
\begin{equation*}
( r^{-k} E_r^l \,  \Phi_{m+1}) (\cdot, \cdot; r) \to 0,
\; \;  \text{as} \; \; r \to 0, 
\quad \text{in the topology of} \; C^{\infty} (F_{p,n}).
\end{equation*} 
\end{proposition}

\begin{proof}
The $N$-th remainder term $e_N (t)$ 
in the Taylor expansion of $e^t$ is expressed as
\begin{equation}\label{smoothness-proposition:2-1}
e_N (t) = t^N g_N (t),
\quad g_N (t) = \frac{1}{(N-1)!} \int_0^1 \, (1-s)^N e^{st} \, ds.
\end{equation}
Thus we have
\begin{equation}
\label{smoothness-proposition:2-2}
\begin{split}
\Phi_{m+1} (\sigma, \omega ; r)
 &= \int_{\sigma \subset \eta \subset \omega^{\perp}} \;
   \int_{y \in \eta^{\perp}} \;
    (-i \langle y, r \omega \rangle )^{m+1}
     g_{m+1} (-i \langle y, r \omega \rangle ) 
       \varphi (\eta, y)
        \, dy d \eta\\
 &= r^{m+1} \, \int_{\sigma \subset \eta \subset \omega^{\perp}} \;
      \int_{y \in \eta^{\perp}} \;
       (-i \langle y,  \omega \rangle )^{m+1}
         g_{m+1} (-i \langle y, r \omega \rangle ) 
          \varphi (\eta, y)
            \, dy d \eta.
\end{split}
\end{equation}
Let us write
\begin{equation}
\label{smoothness-proposition:2-3}
r^{-k} E_r^l \, \left\{
   r^{m+1} (-i \langle y,  \omega \rangle )^{m+1}
       g_{m+1} (-i \langle y, r \omega \rangle ) \right\}
         \varphi (\eta, y)
= r \varphi_{m+1} (\eta, y; \omega, r)
\end{equation}
Here we note that $\varphi_{m+1} (\eta, y; \omega, r)$ 
and its (higher order) derivatives with respect to $y$ are
bounded and integrable in $y$ if $|r| \leq 1$.
By (\ref{smoothness-proposition:2-2}) and (\ref{smoothness-proposition:2-3}), we have
\begin{equation*}
( r^{-k} E_r^l \,  \Phi_{m+1}) (u \cdot \sigma, u \cdot \omega ; r)
= r \int_{\sigma \subset \eta \subset \omega^{\perp}} \;
      \int_{y \in \eta^{\perp}} \;
        \varphi_{m+1} (\eta, u \cdot y; \omega, r)
            \, dy d \eta,
\qquad \text{for} \; \; u \in O(n).
\end{equation*}
Thus we have for $X_1, \cdots , X_L \in so(n)$
\begin{equation}\label{smoothness-proposition:2-5}
(X_1 \cdots X_L \cdot ( r^{-k} E_r^l \, \Phi_{m+1})) (\sigma, \omega ; r)
= r \int_{\sigma \subset \eta \subset \omega^{\perp}} \;
      \int_{y \in \eta^{\perp}} \;
        ((X_1)_y \cdots (X_L)_y \cdot \varphi_{m+1}) (\eta, y; \omega, r)
            \, dy d \eta.
\end{equation}
In equality (\ref{smoothness-proposition:2-5}),
$(X_j)_y, \; (1 \leq j \leq L)$ acts on a function of $y$.
Taking into account that $\varphi \in \mathcal{S} (G(q,n))$
and $g_{m+1} (-it)$ has bounded derivatives,
we see easily that 
\begin{equation}\label{smoothness-proposition:2-6}
\sup_{|r| \leq 1} \, 
  |((X_1)_y \cdots (X_L)_y \cdot \varphi_{m+1}) 
     (\eta, y; \omega, r)|
    \; \leq C (1+||y||)^{- (n-q+1)}, \quad
\text{for} \; y \in \eta^{\perp}.
\end{equation}
In inequality (\ref{smoothness-proposition:2-6}), 
the constant $C$ does not depent on $(\eta, \omega) \in F_{q,n}$.
As a result, for any $L \geq 1$ and for any $X_1, \cdots , X_L \in so(n)$,
there exists a constant $C$ such that
\begin{equation}
\sup_{(\sigma, \omega) \in F_{p,n}} \, 
|(X_1 \cdots X_L \cdot ( r^{-k} E_r^l \,  \Phi_{m+1})) (\sigma, \omega ; r)|
 \leq C |r|, \qquad \text{for} \; \; r , \; (|r|<1),
\end{equation}
which proves the assertion.
\end{proof}

Lemma \ref{T:associatedbundle} and Proposition \ref{smoothness-proposition:2}
yield the following.

\begin{proposition}\label{smoothness-proposition:3}
Let $Q$ be an element of $\u (so(n))$
and let $P =P^{(m)} (E_1, \cdots, E_n)$ be
a homogeneous element of $\u (\mathbb{R}^n)$ of degree $m$.
Then, 
\begin{equation*}
\sup_{(\sigma, \omega) \in F_{p,n}} \, 
  |(QP \cdot \square^{(p)} \Phi_{m+1}) (\sigma, \omega ; r)| \to 0, 
    \quad  \text{as} \; \; r \to 0.
\end{equation*}
\end{proposition}

\begin{proof}
By Theorem \ref{Diff-Op-theorem:1}, the differential operator
$P^{(m)} (E_1, \cdots, E_n)$ on $F_{p,n} \times (\mathbb{R} \setminus \{ 0 \})$
is written in the form
\begin{equation*}
P^{(m)} (E_1, \cdots, E_n) \; = \;
 \sum_{k=0}^m \; A_k^{(m)} (D_{(\sigma, \omega)}) r^{-k} E_r^{m-k} ,
 \quad \text{for} \; \; (\sigma, \omega) \in F_{p,n}, \;
r \neq 0,
\end{equation*}
where $A_k^{(m)} (D_{(\sigma, \omega)})$ is a differential operator
on $F_{p,n}$ of order at most $k$.
Since $Q \in \u (so(n))$, $Q$ can be regarded as a differential operator on $F_{p,n}$.
Moreover, $\square^{(p)}$ acts on functions of $(\sigma, \omega) \in F_{p,n}$ and
therefore $\square^{(p)}$ commutes with the multiplication operator $r$ and
the radial derivative $E_r$.
Taking account of these two facts, we have
\begin{equation*}
(QP \cdot \square^{(p)} \Phi_{m+1}) (\sigma, \omega ; r) =
  \sum_{k=0}^m \; 
    ( \, Q \circ A_k^{(m)} (D_{(\sigma, \omega)}) \circ \square^{(p)}
       \; ( r^{-k} E_r^{m-k} \Phi_{m+1}) \, ) (\sigma, \omega ; r).
\end{equation*}
By Lemma \ref{T:associatedbundle}, $\square^{(p)}$ is a continuous linear operator
on $C^{\infty} (F_{p,n})$ 
and so is $Q \circ A_k^{(m)} (D_{(\sigma, \omega)}) \circ \square^{(p)}$.
Thus, by Proposition \ref{smoothness-proposition:2},
we have
\begin{equation*}
( \, Q \circ A_k^{(m)} (D_{(\sigma, \omega)}) \circ \square^{(p)}
   \; ( r^{-k} E_r^{m-k} \Phi_{m+1}) \, ) (\cdot, \cdot ; r)
     \to 0, \quad \text{as} \; \; r \to 0, \qquad
\text{in the topology of} \; C^{\infty} (F_{p,n}).
\end{equation*}
In particular, we have
\begin{equation*}
\sup_{(\sigma, \omega) \in F_{p,n}} \, 
  |( \, Q \circ A_k^{(m)} (D_{(\sigma, \omega)}) \circ \square^{(p)}
       \; ( r^{-k} E_r^{m-k} \Phi_{m+1}) \, ) (\sigma, \omega ; r) |
         \to 0, \quad \text{as} \; \; r \to 0,
\end{equation*}
which completes the proof.
\end{proof}

\medskip

{\bf Proof of Proposition \ref{T:f-smoothness}. }

\medskip

From now on, we will prove the smoothness of $\widetilde f$
around the set $G_{p,n} \times \{ 0 \}$.
Since $\widetilde f$ is smooth 
on $G(p,n) \setminus (G_{p,n} \times \{ 0 \} )$,
the proof is reduced to prove the following.

\begin{proposition}\label{smoothness-proposition:4}
$\widetilde f$ is smooth 
near $G_{p,n} \times \{ 0 \} $.
\end{proposition}

\begin{proof}
It suffices to show that $\widetilde f$ is of class $C^m$ 
for any $m$ ($m=0,1,2,\cdots,$).
We prove this by induction on $m$.
We have already proved that $\widetilde f$ is continuous on $G(p,n)$,
namely, $\widetilde f$ is of class $C^0$ on $G(p,n)$.
In addition, we have already shown that 
$\widetilde f$ is smooth on $G(p,n) \setminus (G_{p,n} \times \{ 0 \} )$.

We assume that $\widetilde f$ is of class $C^m$
on the set $\{ (\sigma, x) \in G(p,n) \, | \, ||x|| <1 \}$.
Let us take any $U \in \u (\m (n))$ of degree $m+1$.
From now on, we will prove that
$U \cdot {\widetilde f} (\sigma, 0)$ exists and
$U \cdot {\widetilde f} (\sigma, 0)$ is continuous at $(\sigma, 0)$
for any point $(\sigma, 0) \in G_{p,n} \times \{ 0 \}$.

Taking account of Proposition \ref{Diff-Op-lemma:0},
we may assume that $U$ is written in the form (i) or (ii) below.
\begin{itemize}
\item[(i)] \;
$U = X Q P$, where $X \in so(n)$, $Q \in \u (so(n))$ is of degree $k$
and $P \in \u (\mathbb{R}^n)$ is a homogeneous element of degree $m-k$.
\item[(ii)] \;
$U = E_1 P$, where $P \in \u (\mathbb{R}^n)$ is a homogeneous element of degree $m$.
\end{itemize}
So if we prove the following two lemmas, Lemma \ref{smoothness-lemma:1}
and Lemma \ref{smoothness-lemma:2}, then we see that
$\widetilde f$ is of class $C^{m +1}$,
and therefore, the proof of Proposition \ref{smoothness-proposition:4}
is completed.
\end{proof}

\begin{lemma}\label{smoothness-lemma:1}
Let $Q \in \u (so (n))$ be of degree $k$ and let $P \in \u (\mathbb{R}^n)$
be a homogeneous element of degree $l$. We assume that $k+l=m$.
Let $X \in so(n)$ and let $\sigma \in G_{p,n}$.
Then, $Q P \cdot {\widetilde f}$ is differentiable at 
$(\sigma, 0) \in G_{p,n} \times \{ 0 \}$ with respect to $X$.
Moreover, $XQP \cdot {\widetilde f}$ is continuous at $(\sigma, 0)$.
\end{lemma}

\begin{proof}
By the hypothesis of the induction, $Q P \cdot {\widetilde f}$ is
a continuous function on $\{ (\sigma, x) \in G(p,n) \, | \, ||x|| <1 \}$.
By Proposition \ref{smoothness-proposition:1},
\begin{equation*}
{\widetilde f} (\sigma, r \omega ) =
\sum_{j=0}^{l} \; \frac{(-i)^j}{j!} \, P_j (\sigma, x)
  + (\square^{(p)} \, \Phi_{l+1} ) (\sigma, \omega ; r),
  \qquad \text{for} \; \; x = r \omega.
\end{equation*}
Since $P_j (\sigma, x)$ is a homogeneous polynomial of deree $j$ and 
$P \in \u (\mathbb{R}^n)$ is a homogeneous differential operator of order $l$,
$P \cdot P_j (\sigma, x) =0$ if $j \leq l-1$.
Thus we have
\begin{equation}\label{smoothness-lemma:1-2}
QP \cdot {\widetilde f} (\sigma, x)
  = \frac{(-i)^l}{l!} \, (QP \cdot P_l) (\sigma, x)
    + (QP \cdot \square^{(p)} \, \Phi_{l+1} ) (\sigma, \omega ; r),
    \qquad \text{for} \; \; x = r \omega.
\end{equation}
Let $r \to 0$ in (\ref{smoothness-lemma:1-2}).
Then by Proposition \ref{smoothness-proposition:3} and the continuity of 
$Q P \cdot {\widetilde f}$, we have
\begin{equation*}
QP \cdot {\widetilde f} (\sigma, 0)
  = \frac{(-i)^l}{l!} \, (QP \cdot P_l) (\sigma, 0).
\end{equation*}
Hence,
\begin{equation*}
\begin{split}
X \cdot (QP \cdot {\widetilde f}) (\sigma, 0)
 &= \frac{d}{dt} \, (QP \cdot {\widetilde f}) (e^{-tX} \sigma, 0)|_{t=0}\\
 &= \frac{(-i)^l}{l!} \, \frac{d}{dt} \, (QP \cdot P_l) (e^{-tX} \sigma, 0)|_{t=0}\\
 &= \frac{(-i)^l}{l!} \, (XQP \cdot P_l) (\sigma, 0).
\end{split}
\end{equation*}
Note that by the definition of moment condition $P_l \in C^{\infty} (G(p,n))$.
Therefore, $Q P \cdot {\widetilde f}$ is differentiable at 
$(\sigma, 0) \in G_{p,n} \times \{ 0 \}$ with respect to $X$.

Next, we will show the continuity of $XQP \cdot {\widetilde f}$ at $(\sigma, 0)$.
Similarly as in (\ref{smoothness-lemma:1-2}),
\begin{equation*}
XQP \cdot {\widetilde f} (\sigma, x)
  = \frac{(-i)^l}{l!} \, (XQP \cdot P_l) (\sigma, x)
    + (XQP \cdot \square^{(p)} \, \Phi_{l+1} ) (\sigma, \omega ; r),
      \qquad \text{for} \; \; x = r \omega.
\end{equation*}
Again, by Proposition \ref{smoothness-proposition:3},
\begin{equation*}
\begin{split}
\lim_{x \to 0} \, (XQP \cdot {\widetilde f}) (\sigma, x)
 &= \frac{(-i)^l}{l!} \, \lim_{x \to 0} \, (XQP \cdot P_l) (\sigma, x)
    + \lim_{x \to 0} \, (XQP \cdot \square^{(p)} \, \Phi_{l+1} ) (\sigma, \omega ; r)\\
 &= \frac{(-i)^l}{l!} \, (XQP \cdot P_l) (\sigma, 0) 
     = (XQP \cdot {\widetilde f}) (\sigma, 0).
\end{split}
\end{equation*}
\end{proof}

\begin{lemma}\label{smoothness-lemma:2}
Let $P \in \u (\mathbb{R}^n)$
be a homogeneous element of degree $m$.
Then, $ P \cdot {\widetilde f}$ is differentiable at 
$(\sigma, 0) \in G_{p,n} \times \{ 0 \}$ with respect to $E_1$.
Moreover, $E_1 P \cdot {\widetilde f}$ is continuous at $(\sigma, 0)$.
\end{lemma}

\begin{proof}
If $\mathbf{e}_1$ is parallel to $\sigma$, obviously
the assertion holds. so we assume that $\mathbf{e}_1 \nparallel \sigma$.
By the hypothesis of the induction, $P \cdot {\widetilde f}$ is continuous
on $\{ (\sigma, x) \in G(p,n) \, | \, ||x|| <1 \}$.
Let $v=\Pr_{\sigma^{\perp}} \mathbf{e}_1 \neq 0$.
By the mean value theorem,
\begin{equation}\label{smoothness-lemma:2-1}
\frac{1}{t} 
\{ P \cdot {\widetilde f} (\sigma, -t v)
    - P \cdot {\widetilde f} (\sigma, 0) \}
= (E_v P \cdot {\widetilde f}) (\sigma, -\theta t v)
= (E_1 P \cdot {\widetilde f}) (\sigma, -\theta t v),
\end{equation}
for some $\theta, \; (0 < \theta < 1)$.

By Proposition \ref{smoothness-proposition:1},
\begin{equation}\label{smoothness-lemma:2-2}
{\widetilde f} (\sigma, r \omega ) =
\sum_{j=0}^{m+1} \; \frac{(-i)^j}{j!} \, P_j (\sigma, x)
  + (\square^{(p)} \, \Phi_{m+2} ) (\sigma, \omega ; r),
  \qquad \text{for} \; \; x = r \omega.
\end{equation}
Since $P_j (\sigma, x)$ is a homogeneous polynomial of degree $j$ and 
$E_1 P \in \u (\mathbb{R}^n)$ is a homogeneous differential operator of order $m+1$,
$E_1 P \cdot P_j (\sigma, x) =0$ if $j \leq m$.
Thus we have
\begin{equation}\label{smoothness-lemma:2-3}
(E_1 P \cdot {\widetilde f}) (\sigma, r \omega ) =
 \frac{(-i)^{m+1}}{(m+1)!} \, (E_1 P \cdot P_{m+1}) (\sigma, x)
  + ( E_1 P \circ \square^{(p)} \, \Phi_{m+2} ) (\sigma, \omega ; r),
  \qquad \text{for} \; \; x = r \omega.
\end{equation}
Hence by Proposition \ref{smoothness-proposition:3},
\begin{equation}\label{smoothness-lemma:2-4}
\begin{split}
(E_1 P \cdot {\widetilde f}) (\sigma, -\theta t v)
 &= \frac{(-i)^{m+1}}{(m+1)!} \, (E_1 P \cdot P_{m+1}) (\sigma, -\theta t v)\\
 & \; \; + ( E_1 P \circ \square^{(p)} \, \Phi_{m+2} ) (\sigma, v/||v|| ; - \theta t ||v||)\\
 &\longrightarrow \frac{(-i)^{m+1}}{(m+1)!} \, (E_1 P \cdot P_{m+1}) (\sigma, 0),
   \qquad \text{as} \; \; t \to 0.
\end{split}
\end{equation}
Therefore, it follows 
from (\ref{smoothness-lemma:2-1}) and (\ref{smoothness-lemma:2-4}) that 
$ P \cdot {\widetilde f}$ is differentiable at 
$(\sigma, 0) \in G_{p,n} \times \{ 0 \}$ with respect to $E_1$
and that
\begin{equation}\label{smoothness-lemma:2-5}
E_1 P \cdot {\widetilde f} (\sigma, 0)
= \frac{(-i)^{m+1}}{(m+1)!} \, (P \cdot P_{m+1}) (\sigma, 0).
\end{equation}

Next we will show that $(E_1 P \cdot {\widetilde f})$
is continuous at $(\sigma, 0)$.
By (\ref{smoothness-lemma:2-3}) 
and Proposition \ref{smoothness-proposition:3},
\begin{equation}\label{smoothness-lemma:2-6}
\begin{split}
(E_1 P \cdot {\widetilde f}) (\sigma, r \omega ) 
 &= \frac{(-i)^{m+1}}{(m+1)!} \, (E_1 P \cdot P_{m+1}) (\sigma, x)\\
 & \; \; + ( E_1 P \circ \square^{(p)} \, \Phi_{m+2} ) (\sigma, \omega ; r)\\
 & \longrightarrow \frac{(-i)^{m+1}}{(m+1)!} \, (E_1 P \cdot P_{m+1}) (\sigma, 0)
  = (E_1 P \cdot {\widetilde f}) (\sigma, 0 ),
  \qquad \text{as} \; \; r \to 0,
\end{split}
\end{equation}
which proves the continuity of $(E_1 P \cdot {\widetilde f})$
at $(\sigma, 0)$.
\end{proof}

Finally we will prove that ${\widetilde f}$ is a Schwartz class function
on $G(p,n)$.

Before going to the proof, let us recall
from \cite{Gonzalez-Kakehi-1} and \cite{Richter}
that a smooth function $g$
on $G(p,n)$ belongs to the Schwartz space
${\mathcal S} (G(p,n))$ if for any nonnegative integers $N$ and $m$ and 
for any $m$ vector fields $Y_1, \cdots, Y_m \; \in \m (n)$ g satisfies
$$
\sup_{(\sigma, x) \in G(p,n)} \; 
  ||x||^N \, |( Y_1 \cdots Y_m \cdot g) (\sigma, x)|
\; < \; \infty.
$$

\medskip

{ \bf Proof of Proposition \ref{T:f-rapid-decrease}.}

\begin{proof}
We start with the following.
\begin{equation*}
{\widetilde f} (\sigma, r \omega) 
= \square^{(p)} \int_{\sigma \subset \eta \subset \omega^{\perp}} \;
      \mathcal{F}_q \varphi (\eta, r \omega) \, d \eta.
\end{equation*}
Let us take any nonnegative integers $m$ and $N$.
In addition, let us take
any element $Q \in \u (so(n))$
and any homogeneous element
$P =P^{(m)} (E_1, \cdots, E_n) \in \u (\mathbb{R}^n)$ of degree $m$.
By Theorem \ref{Diff-Op-theorem:1},
\begin{equation}\label{smoothness-proposition:5-2}
\begin{split}
{}& r^N \, (Q P^{(m)} (E_1, \cdots, E_n) \cdot {\widetilde f}) (\sigma, r \omega) \\
&= r^N \, \sum_{k=0}^m \; 
   Q A_k^{(m)} (D_{(\sigma, \omega)}) r^{-k} E_r^{m-k} \, 
    \square^{(p)} \, \int_{\sigma \subset \eta \subset \omega^{\perp}} \;
      \mathcal{F}_q \varphi (\eta, r \omega) \, d \eta \\
&= \sum_{k=0}^m \; 
  Q A_k^{(m)} (D_{(\sigma, \omega)}) 
    \square^{(p)} \, \int_{\sigma \subset \eta \subset \omega^{\perp}} \;
      r^{N-k} E_r^{m-k} \, \mathcal{F}_q \varphi (\eta, r \omega) \, 
        d \eta \\
&= \sum_{k=0}^m \; 
   ( \, Q \circ A_k^{(m)} (D_{(\sigma, \omega)}) 
    \circ \square^{(p)} \circ S \, ) 
      ( \, r^{N-k} E_r^{m-k} \, \mathcal{F}_q \varphi (\cdot , r \cdot ) \, )
        (\sigma, \omega).
\end{split}
\end{equation}
where $A_k^{(m)} (D_{(\sigma, \omega)})$ is a differential operator
on $F_{p,n}$ of order at most $k$ and
where $S$ is the Radon transform from $C^{\infty} (F_{q,n})$ to
$C^{\infty} (F_{p,n})$ defined by (\ref{E:flagint2}).
In the above computation, we used the fact that
the multiplication operators $r^j \; (j = 1,2, \cdots)$ commute with
differential operators on $F_{p,n}$ and with the operator $\square^{(p)}$.
Since $\mathcal{F}_q \varphi \in \mathcal{S} (G(q,n))$,
\begin{equation*}
r^{N-k} E_r^{m-k} \, \mathcal{F}_q \varphi (\cdot , r \cdot )
\longrightarrow 0, \qquad \text{as} \; \; |r | \to \infty,
\end{equation*}
in the topology of $C^{\infty} (F_{q,n})$.
Moreover, in the summand of (\ref{smoothness-proposition:5-2}),
$Q \circ A_k^{(m)} (D_{(\sigma, \omega)}) \circ \square^{(p)} \circ S$
is a continuous linear operator from
$C^{\infty} (F_{q,n})$ to $C^{\infty} (F_{p,n})$.
Thus we have
\begin{equation*}
( \, Q \circ A_k^{(m)} (D_{(\sigma, \omega)}) 
    \circ \square^{(p)} \circ S \, ) 
      ( \, r^{N-k} E_r^{m-k} \, \mathcal{F}_q \varphi (\cdot , r \cdot ) \, )
\longrightarrow 0, \qquad \text{as} \; \; |r | \to \infty,
\end{equation*}
in the topology of $C^{\infty} (F_{p,n})$,
from which we can conclude that
\begin{equation}\label{smoothness-proposition:5-5}
\sup_{(\sigma, \omega) \in F_{p,n}} \;
| r^N \, (Q P^{(m)} (E_1, \cdots, E_n) \cdot {\widetilde f}) (\sigma, r \omega) |
\longrightarrow 0, \qquad \text{as} \; \; |r| \to \infty.
\end{equation}
Taking account of Proposition \ref{Diff-Op-lemma:0}, we see that
(\ref{smoothness-proposition:5-5}) proves
the assertion.
\end{proof}

\end{section}


\begin{section}{The Support Theorem I}

Our objective in this section is to prove a support theorem for
$\mathcal R^{(p,q)}$ based on the forward moment conditions ($H'$), generalizing
Theorem~\ref{T:classicalsupport}.  For this, we make use of the
second-order differential operator $\Delta_p$ on $G(p,n)$ which acts
as the Laplacian on each fiber $\sigma^\perp$:
$$
\Delta_p f(\sigma,x)\;=\;L_{\sigma^\perp} f(\sigma,x)
$$
$\Delta_p$ is invariant under the motion group $M(n)$
(\cite{Helgason-1965}).  We define the differential operator
$\Delta_q$ on $G(q,n)$ similarly.

Using the operator $\Delta_p$, it makes sense to talk about harmonic
polynomials on the fiber $\sigma^\perp$ in $G(p,n)$, as well as
spherical harmonics on the unit sphere $S_\sigma$ in $\sigma^\perp$.
(These can also be obtained from harmonic polynomials on $\mathbb
R^{n-p}$ by a translation by an appropriate element $u\in O(n)$.)

For any $R>0$, define $\bar\beta_p^R$ to be the set of all $p$-planes
$(\sigma,x)$ at distance $\leq R$ from the origin.  We will use the
following analogue of Theorem~\ref{T:Paleypolar1} for the partial
Fourier transform $\mathcal F_p$.

\begin{theorem}\label{T:Paleypolar2}
The partial Fourier transform
$f\mapsto \mathcal F_p f$ maps $\mathcal D(\bar\beta_p^R)$ onto the
space of all functions $\mathcal F_p
f(\sigma,\lambda\omega)=\widetilde F(\sigma,\omega;\lambda)\in
C^\infty(F_{p,n}\times\mathbb R)$ satisfying the following conditions:
\item{(i)} For each $(\sigma,\omega)\in F_{p,n}$, the function
$\lambda\mapsto \widetilde F(\sigma,\omega;\lambda)$ extends to a holomorphic
function on $\mathbb C$ with the property that
\begin{equation}\label{E:uniform-paley-estimate}
\sup_{(\sigma,\omega;\lambda)\in F_{p,n}\times\mathbb C}
\left|(1+|\lambda|)^N\,\widetilde F(\sigma,\omega;\lambda)
 \, e^{-R|\text{Im}\lambda|}\right|\;<\;\infty,
\end{equation}
for each $N\in\mathbb Z^+$.
\item{(ii)} For each $k\in\mathbb Z^+$, for each $\sigma\in G_{p,n}$,
and for each homogeneous degree $k$ harmonic polynomial $h_\sigma$ on
$\sigma^\perp$, the function
$$
\lambda\;\mapsto\;\lambda^{-k}\,\int_{S_\sigma}
\widetilde F(\sigma,\omega;\lambda) h_\sigma(\omega)\,d\omega
$$
is even and holomorphic in $\mathbb C$.
\end{theorem}

\begin{proof}
If $f\in\mathcal D(\bar\beta_p^R)$, then clearly the function
$\lambda\mapsto\widetilde F(\sigma,\omega;\lambda)=
\int_{\sigma^\perp} f(\sigma,x) e^{-i\lambda\langle
x,\omega\rangle}\,dx$ extends to a holomorphic function on $\mathbb
C$.  If $N\in\mathbb Z^+$, we have
$$
\begin{aligned}
| \lambda |^{2N} \,
\left| \widetilde F(\sigma,\omega;\lambda)\right|
\;&=\;\left|\lambda^{2N} \,
  \int_{\sigma^\perp} f(\sigma,x)\,e^{-i\lambda\langle
    x,\omega\rangle}\,dx \right|\\
&=\;\left|\int_{\sigma^\perp} (-\Delta_p)^N f(\sigma,x)\,
  e^{-i\lambda\langle x,\omega\rangle}\,dx \right|\\
&\leq\;  C
  \max_{(\sigma, x) \in \bar\beta_p^R} \, |(-\Delta_p)^N f (\sigma, x)| 
    \,e^{R|\text{Im}\lambda|} 
\end{aligned}
$$
for all $(\sigma,\omega;\lambda)\in F_{p,n}\times \mathbb C$, which
proves the estimate (\ref{E:uniform-paley-estimate}).  Moreover, $\widetilde F$ satisfies
condition (ii) above by Theorem~\ref{T:Paleypolar1} applied to the Euclidean
space $\sigma^\perp$.

Conversely, suppose $\mathcal F_p f\,=\,\widetilde F$ satisfies (i) and (ii).
Let $\sigma\in G_{p,n}$.  Then -- again by Theorem~\ref{T:Paleypolar1} applied to the Euclidean
space $\sigma^\perp$ -- the function
$$
f(\sigma,x)\;=\;(2\pi)^{p-n}\int_{S_\sigma}\,\int_0^\infty \widetilde 
F(\sigma,\omega;\lambda)\,e^{i\lambda\langle
x,\omega\rangle}\,\lambda^{n-p-1}\,d\lambda\,d\omega
$$
satisfies $f(\sigma,x)=0$ for all $||x||>R$.  This shows that $f$ is supported in
$\bar\beta_p^R$.
\end{proof}

Our main result, the support theorem below, extends the classical
support theorem, Theorem~\ref{T:classicalsupport}.

\begin{theorem}\label{T:Grassmanniansupport} (Support Theorem for
$\mathcal R^{(p,q)}$.)
Assume that $\rank(G(p,n))\,=\,\rank(G(q,n))$.  Suppose that $f\in\mathcal
S(G(p,n))$ satisfies $\mathcal R^{(p,q)} f(\eta,y)=0$ whenever
$|y|>R$.  Then $f$ is supported in $\bar\beta_p^R$.
\end{theorem}

\begin{proof}
Let $\varphi=\mathcal R^{(p,q)} f$; then let
$\widetilde\Phi(\eta,\omega;\lambda)=
{\mathcal F}_q \varphi(\eta,\lambda\omega)$ as in (\ref{E:definePhi}).
From Theorem~\ref{T:Paleypolar2}, the function
$\lambda\mapsto\widetilde\Phi(\eta,\omega;\lambda)$ extends to a
holomorphic function on $\mathbb C$ for each $(\eta,\omega)\in
F_{q,n}$, and satisfies the estimate
\begin{equation}\label{E:Phi-estimate}
\sup_{(\eta,\omega;\lambda)\in F_{q,n}\times\mathbb C}
\left|(1+|\lambda|)^N\,\widetilde\Phi(\eta,\omega;\lambda)
\, e^{-R|\text{Im}\lambda|}\right| \; < \; \infty,
\qquad \text{for all} \; \; N\in\mathbb Z^+.
\end{equation}
Moreover, since 
$X \cdot {\mathcal F}_q \varphi = {\mathcal F}_q (X \cdot \varphi)$
for $X \in so(n)$,
we have a similar estimate for $X_1 \cdots X_m \cdot \widetilde\Phi$.
Here $X_1, \cdots, X_m \in so(n)$. 
Namely, for any nonnegative integer $m$ and for any $m$
vector fields $X_1, \cdots, X_m \in so(n)$, we have
\begin{equation}\label{E:support-theorem-2-3:PW-estimate}
\sup_{(\eta,\omega;\lambda)\in F_{q,n}\times\mathbb C}
\left|(1+|\lambda|)^N\,
(X_1 \cdots X_m \cdot \widetilde\Phi) (\eta,\omega;\lambda)
\, e^{-R|\text{Im}\lambda|}\right| \; < \; \infty,
\qquad \text{for all} \; \; N\in\mathbb Z^+.
\end{equation}
Let us fix an arbitrary nonnegative integer $N$ and introduce
a family of functions 
$\{ H_{\lambda} \; | \; \lambda \in \mathbb{C} \, \}$
in $C^{\infty} (F_{q,n})$ as follows.
\begin{equation}
H_{\lambda} (\eta,\omega) = 
  \widetilde\Phi (\eta,\omega;\lambda) (1+|\lambda|)^N e^{-R|\text{Im}\lambda|}.
\end{equation}
Then the above two estimates 
(\ref{E:Phi-estimate}) and (\ref{E:support-theorem-2-3:PW-estimate})
show that $\{ H_{\lambda} \; | \; \lambda \in \mathbb{C} \, \}$
is a bounded set in $C^{\infty} (F_{q,n})$.

Now we put $\widetilde F(\sigma,\omega;r)\,=\,\mathcal F_p f(\sigma,r\omega)$
as in (\ref{E:define-f}).
Then by the (projection-slice) Theorem~\ref{E:projectionslice}, 
$$
\widetilde\Phi(\eta,\omega;r)\;=\;\int_{\sigma\subset\eta} \widetilde
F(\sigma,\omega;r)\,d\sigma.
$$
(This is the same as equation~(\ref{E:flagint1}).) 
Just as with
(\ref{E:widetildeF}), for fixed $\omega$ 
the inversion formula (\ref{E:Rpqomegainversion}) may again be applied to
recover $\widetilde F$ from $\widetilde\Phi$:
$$
\widetilde F(\sigma,\omega;r)\;=\;\square_{p,q}^\omega(
S\widetilde\Phi)(\sigma,\omega;r)\;
  =\;\square^{(p)} \circ S \widetilde\Phi (\sigma,\omega;r).
$$
This formula still holds when $r$ is replaced by a complex parameter
$\lambda$.
\begin{equation}\label{E:widetildeF2}
\widetilde F(\sigma,\omega;\lambda) \; 
  =\; \square^{(p)} \circ S \widetilde\Phi (\sigma,\omega;\lambda).
\end{equation}
Next, we put
$G_{\lambda} (\sigma,\omega) 
 = \widetilde F(\sigma,\omega;\lambda) (1+|\lambda|)^N e^{-R|\text{Im}\lambda|}$.
Then $G_{\lambda} = \square^{(p)} \circ S H_{\lambda}$.
Since by Lemma \ref{T:S-smoothness} and Lemma \ref{T:associatedbundle}
$\square^{(p)} \circ S$ is a continuous linear operator
from $C^{\infty} (F_{q,n})$ to $C^{\infty} (F_{p,n})$,
the set $\{ G_{\lambda} \; | \; \lambda \in \mathbb C \}$ is bounded
in $C^{\infty} (F_{p,n})$.
In particular, we have
\begin{equation*}
 \sup_{(\sigma,\omega;\lambda)\in F_{p,n}\times\mathbb C} \; 
  | G_{\lambda} (\sigma,\omega) |
= \sup_{(\sigma,\omega;\lambda)\in F_{p,n}\times\mathbb C} \; 
  | \widetilde F(\sigma,\omega;\lambda) (1+|\lambda|)^N e^{-R|\text{Im}\lambda|} |
  \; < + \infty,
\end{equation*}
which shows that the function 
$\lambda \mapsto \widetilde F(\sigma,\omega;\lambda)$
satisfies the uniform Paley-Wiener estimate (\ref{E:uniform-paley-estimate}).

Next, we will prove that
the function $\lambda \mapsto \widetilde
F(\sigma,\omega;\lambda)$ is holomorphic in $\mathbb C$ for each
$(\sigma,\omega)\in F_{p,n}$.
Let us take any closed curve $\Gamma$ in $\mathbb C$.
Then by the continuity of $\square^{(p)} \circ S$,
$$
\int_{\Gamma} \; F(\sigma,\omega;\lambda) \, d \lambda
= \square^{(p)} \circ S \left(
   \int_{\Gamma} \widetilde\Phi (\cdot,\cdot;\lambda) d \lambda
     \right) (\sigma,\omega) = 0.
$$
Here we used the fact that $\widetilde\Phi$ is holomorphic in
$\lambda$.
Therefore, by Morera's theorem, 
$\lambda \mapsto \widetilde F (\sigma,\omega;\lambda)$ is holomorphic.

To finish the proof of Theorem~\ref{T:Grassmanniansupport} it remains to
prove that for each homogeneous spherical harmonic $h_\sigma$ of
degree $k$ on the unit sphere $S_\sigma$ in $\sigma^\perp$, the
function
\begin{equation}\label{E:harmonic-int}
\lambda\;\mapsto\;\lambda^{-k}\,\int_{S_\sigma}\widetilde
F(\sigma,\omega;\lambda)\,
h_\sigma(\omega)\,d\omega
\end{equation}
is even  and holomorphic in $\lambda$.  First, the fact that $\widetilde
F(\sigma,\omega;\lambda)=\widetilde F(\sigma,-\omega,-\lambda))$ implies
that (\ref{E:harmonic-int}) is even in $\lambda$.

In the calculations leading to Proposition \ref{smoothness-proposition:1},
which use the forward moment condition ($H'$),
the real parameter $r$ can be replaced by the complex parameter
$\lambda$ to give us the following expression.
$$
\widetilde\Phi(\eta,\omega;\lambda)
= \; \sum_{l=0}^{k-1} 
   \frac{(-i\lambda)^l}{l!}\,\int_{\sigma\subset\eta}
    P_l(\sigma,\omega)\,d\sigma\;
 + \; \int_{\eta^\perp}  \varphi(\eta,y) \,(-i \lambda \langle y,\omega \rangle)^k
    \, g_k(-i\lambda \langle y, \omega \rangle) \, dy,
$$
where $g_k(-iz)\,=\, e_k(-iz)/(-iz)^k$ is holomorphic in $z$ and bounded if
$z$ is real.  
(See (\ref{smoothness-proposition:2-1}) for an explicit expression of $g_k$.)
We write the last integral above as
$$
\lambda^k\,\int_{\eta^\perp}\varphi(\eta,y)\,(-i\langle
y,\omega\rangle)^k
\, g_k(-i\lambda\langle y,\omega\rangle)\,dy\;=\;\lambda^k\,\Psi_k(\eta,\omega;\lambda),
$$
where $\Psi_k(\eta,\omega;\lambda)$ is a smooth function on
$F_{q,n}\times\mathbb C$ and holomorphic in $\lambda$ for each
$(\eta,\omega)\in F_{q,n}$.

By the inversion formula (\ref{E:widetildeF2}) (and
(\ref{E:compactinversion}) for the function $P_l \in C^{\infty} (G (p,n))$),
we get
$$
\widetilde F(\sigma,\omega;\lambda)\;
= \;\sum_{l=0}^{k-1}
   \frac{(-i\lambda)^l}{l!}\,P_l(\sigma,\omega)\;
   + \; \lambda^k\,\square^{(p)} S \Psi_k(\sigma,\omega;\lambda).
$$
Similarly as in the proof of holomorphicity of 
$\widetilde F (\sigma,\omega;\lambda)$,
we can prove that $\square^{(p)}
S\Psi_k(\sigma,\omega;\lambda)$
is holomorphic in $\lambda$ for each $(\sigma,\omega)$.  Hence
\begin{align*}
\int_{S_\sigma} \widetilde
F(\sigma,\omega;\lambda)\,h_\sigma(\omega)\,d\omega 
\;&=\;\sum_{l=0}^{k-1} 
 \frac{(-i\lambda)^l}{l!}\,\int_{S_\sigma}
   P_l(\sigma,\omega)\,h_\sigma(\omega)\,d\omega\;
 + \; \lambda^k\,\int_{S_\sigma}\square^{(p)} S
   \Psi_k(\sigma,\omega;\lambda) \,h_\sigma(\omega) \, d\omega\\
&=\;\lambda^k\,\int_{S_\sigma}\square^{(p)} S
 \Psi_k(\sigma,\omega;\lambda) \, h_\sigma(\omega) \, d\omega,
\end{align*}
since $P_l(\sigma,\omega)$ is a degree $l$ polynomial in $\omega\in
S_\sigma$, and so is a sum of spherical harmonics in $S_\sigma$ of
degree $\leq l$.  The fact that the mapping
$$
\lambda\;\mapsto\;\lambda^{-k}\,\int_{S_\sigma} \widetilde
F(\sigma,\omega;\lambda)\,h_\sigma(\omega)\,d\omega\;=\;\int_{S_\sigma}
\square^{(p)} S \Psi_k(\sigma,\omega;\lambda)\,h_\sigma(\omega)\,d\omega
$$
is holomorphic in $\lambda$ now follows.  The function 
$\mathcal F_p f(\sigma,\lambda\omega)\,=\,\widetilde
F(\sigma,\omega;\lambda)$ thus satisfies conditions (i) and (ii) in 
Theorem~\ref{T:Paleypolar2}, and hence by that theorem $f$ is supported in
$\bar\beta_p^R$.  This completes the
proof of Theorem~\ref{T:Grassmanniansupport}.
\end{proof}

\end{section}


\begin{section}{Support Theorem II.}

In this section, we prove the support theorem for the Radon transform
$\mathcal{R}^{\, (p,q)}: \mathcal{S} (G(p,n)) \to \mathcal{S} (G(q,n))$
in the case when $p <q$ and $\dim G(p,n) < \dim G(q,n)$.

Let $\mathcal {O}$ be a subset in $\mathbb{R}^n$.
Throughout this section, we assume the following condition $(A)$
on $\mathcal {O}$.
\begin{equation}
\text{For any $\ell \in G(p,n)$ with $\ell \subset \mathcal {O}^c$,
there exists a hyperplane $L$
 such that $\ell \subset L \subset \mathcal {O}^c$}.
\tag{A}
\end{equation}

The following support theorem holds.

\begin{theorem}\label{th:support-theorem-1}
Suppose that $f \in \mathcal{S} (G(p,n))$.
If $\mathcal{R}^{\, (p,q)} \, f (\xi) = 0$ for all $\xi \subset \mathcal {O}^c$,
then $f(\ell) = 0$ for all $\ell \subset \mathcal {O}^c$.
\end{theorem}

As a corollary, we have the usual support theorem for $\mathcal{R}^{\, (p,q)}$, namely,

\begin{cor}
Let $r>0$.
Suppose that $f \in \mathcal{S} (G(p,n))$.
If $\mathcal{R}^{\, (p,q)} \, f (\xi) = 0$ for all $\xi, \; \; \rm{dist} ( \xi, 0) > r$,
then $f(\ell) = 0$ for all $\ell, \; \; \rm{dist} (\ell, 0) > r $.
\end{cor}

\begin{proof}
If $\mathcal {O} = \{ \; x \in \mathbb{R}^n \; | \; || x || <r \; \}$,
then obviously $\mathcal {O}$ satisfies the condition $(A)$.
\end{proof}

\begin{remark}
Similarly, if $\mathcal {O}$ is a convex set, then $\mathcal {O}$
satisfies the condition $(A)$.
So, in this case, the support theorem also holds.
\end{remark}

However, there are many cases when $\mathcal {O}$ satisfies the condition $(A)$ but
$\mathcal {O}$ is not necessarily convex.
Even in such a case, the support theorem holds.

The following are examples.

\medskip

{\it Example 1.}
\quad
For $a, b \; (a<b)$ and for $r>0$, let
\begin{equation}
\begin{split}
\mathcal {O}_1 &= 
  \{ x = (x_1, x_2, \cdots, x_n) \in \mathbb{R}^n \; | \;
      (x_1 - a)^2 + x_2^2 + \cdots + x_n^2 \leq r^2, \; x_1 \leq a \; \},\\
\mathcal {O}_2 &= 
  \{ x = (x_1, x_2, \cdots, x_n) \in \mathbb{R}^n \; | \;
      (x_1 - b)^2 + x_2^2 + \cdots + x_n^2 \leq r^2, \; x_1 \geq b \; \},
\end{split}
\end{equation}
and let $\mathcal {O} = \mathcal {O}_1 \cup \mathcal {O}_2$.
Then, $\mathcal {O}$  is no longer convex. But it is easily seen that
$\mathcal {O}$ satisfies the condition $(A)$.

\medskip

{\it Example 2.}
\quad
For $a, b \; (a<b)$, let
\begin{equation}
\mathcal {O}_1 = 
  \{ x = (x_1, x_2, \cdots, x_n) \in \mathbb{R}^n \; | \;
      x_1 \leq a \; \}, \quad
\mathcal {O}_2 = 
  \{ x = (x_1, x_2, \cdots, x_n) \in \mathbb{R}^n \; | \;
      x_1 \geq b \; \},
\end{equation}
and let $\mathcal {O} = \mathcal {O}_1 \cup \mathcal {O}_2$.
Then, $\mathcal {O}^c$ is a band domain. Similarly as in the above example,
$\mathcal {O}$ is not convex, but $\mathcal {O}$ satisfies the condition $(A)$.

Now, we proceed to prove Theorem \ref{th:support-theorem-1}.
The key is the injectivity of the Radon transform.

\medskip

{\bf The proof of Theorem \ref{th:support-theorem-1}.}
\par
Let us take an arbitrary $p$-plane $\ell_0$ which is included in $\mathcal {O}^c$.
Then, by the condition $(A)$, there exists a hyperplane $L \in G(n-1, n)$
such that $\ell_0 \subset L \subset \mathcal {O}^c$.
Let $G(d, n ; L) = \{ \; \gamma \in G(d,n) \; | \; \gamma \subset L \; \}$.
Then the space of Schwartz class functions $\mathcal{S} (G(d,n; L))$
on $G(d, n ; L)$ is defined in a similar manner to $\mathcal{S} (G(d,n))$.
In addition, we can define a Radon transform
$\mathcal{R}_L^{\, (p,q)}: \mathcal{S} (G(p,n; L)) \to \mathcal{S} (G(q,n; L))$
as follows.
\begin{equation}\label{support-theorem-1:1}
\mathcal{R}_L^{\, (p,q)} f (\xi) 
  = \int_{\ell \subset \xi} \;
     f(\ell) \, d_{\xi} \ell,
\quad \text{for} \; \xi \in G(q,n; L),
\quad \text{and for} \; f \in \mathcal{S} (G(q,n; L)),
\end{equation}
where $d_{\xi} \ell$ is the canonical measure on the set
$\{ \ell \in G(p,n ; L) \; | \; \ell \subset \xi \; \}$.
For $f \in \mathcal{S} (G(p,n))$ in the statement of Theorem
\ref{th:support-theorem-1}, let $f_L$ be the restriction of $f$ onto
the submanifold $G(p,n ; L)$. Then, $f_L \in \mathcal{S} (G(p,n; L))$.
Moreover, by the definition of $\mathcal{R}_L^{\, (p,q)}$ and
by the assumption of the Theorem, we have
\begin{equation}\label{support-theorem-1:2}
\mathcal{R}_L^{\, (p,q)} \, f_L (\xi) = \mathcal{R}^{\, (p,q)} \, f (\xi) = 0,
\quad \text{for} \; \xi \in G(q, n; L).
\end{equation}
In fact, if $\xi \in G(q, n; L)$, then $\xi \subset \mathcal{O}^c$ and therefore
$\mathcal{R}^{\, (p,q)} \, f (\xi) = 0$.

By applying a suitable translation and a suitable orthogonal transformation to $L$,
we may assume that 
$L = \mathbb{R} \mathbf{e}_1  \oplus \cdots \oplus \mathbb{R} \mathbf{e}_{n-1}
 (\cong \mathbb{R}^{n-1})$.
Then, the Radon transform $\mathcal{R}_L^{\, (p,q)}$ is nothing but 
the Radon transform from $\mathcal{S} (G(p,n-1))$ to $\mathcal{S} (G(q,n-1))$
associated with inclusion incidense relation.
By the assumption that $p <q$ and $\dim G(p,n) < \dim G(q,n)$,
we see easily that $\dim G(p,n-1) \leq \dim G(q,n-1)$.
So, by Theorem 6.4 and Remark 3 in Section 6 of our previous paper
\cite{Gonzalez-Kakehi-1}, $\mathcal{R}_L^{\, (p,q)}$ is injective.
Hence, by (\ref{support-theorem-1:2}), 
$f_L (\ell) =0$ for $\ell \in G(p,n; L)$.
In particular, $f (\ell_0) = f_L (\ell_0) =0$, which completes the proof.
\qed

Finally, as an application of the above support theorem, we
give a range characterization of $\mathcal{R}^{\, (p,q)}$
in the category of compactly supported smooth functions.

Let $C_c^{\infty} (G(d,n))$ denote the space of compactly supported
smooth functions on $G(d,n)$.
Then, it is easily seen that since $p<q$ the Radon transform $\mathcal{R}^{\, (p,q)}$
maps $C_c^{\infty} (G(p,n))$ to $C_c^{\infty} (G(q,n))$.
As we stated in the introduction, since $\dim G(p,n) < \dim G(q,n)$
the image of $\mathcal{R}^{\, (p,q)}$ is characterized as
the solution space of
an $M(n)$-invariant differential equation of order $2p+4$ of the form,
\begin{equation}
d \nu (Q_{2p+4}) \, \varphi \; = \; 0.
\end{equation}
Here $Q_{2p+4}$ is an element in $\frak{z} (\frak{m} (n))$
and is expressed as the sum of the squares of ^^ ^^ Pfaffians" of order $p+2$.
(See Gonzalez and Kakehi \cite{Gonzalez-Kakehi-1} for
the definition of the operator $Q_{2p+4}$ and its detailed properties.)

Since $C_c^{\infty} (G(d,n)) \subset \mathcal{S} (G(p,n))$,
the support theorem (Theorem \ref{th:support-theorem-1}) and
the range theorem for $\mathcal{R}^{\, (p,q)}$
(Theorem 7.7 of \cite{Gonzalez-Kakehi-1})
yield the following.

\begin{theorem}
Suppose that $p <q$ and $\dim G(p,n) < \dim G(q,n)$.
A function $\varphi \in C_c^{\infty} (G(q,n))$
belongs to the range 
$\mathcal{R}^{\, (p,q)} (C_c^{\infty} (G(p,n)))$
if and only if 
$d \nu (Q_{2p+4}) \, \varphi \; = \; 0$.
\end{theorem}

\end{section}


\end{document}